\documentclass[11pt]{amsart}
\usepackage[usenames,dvipsnames]{xcolor}
\usepackage{appendix}
 \usepackage{amsthm}
 \usepackage{amsmath}
 \usepackage{amsfonts}
 \usepackage{amssymb}
 \usepackage{amscd}
\usepackage{graphicx, subfigure}
\usepackage{mathtools}
\usepackage{hyperref}

\usepackage[normalem]{ulem}

\usepackage{colortbl} 

\usepackage{url}
\usepackage{algorithm}
\usepackage{algorithmic}

\usepackage{float}	

\newcommand{\vcenteredinclude}[2]{\begingroup
\setbox0=\hbox{\includegraphics[height=#2]{#1}}%
\parbox{\wd0}{\box0}\endgroup}

\newtheorem{thm}{Theorem}[section]

\newtheorem{lem}[thm]{Lemma}
\newtheorem{prop}[thm]{Proposition}
\theoremstyle{definition}
\newtheorem{ex}[thm]{Example}
\newtheorem{defn}[thm]{Definition}

\numberwithin{equation}{section}

\def \N { {\mathbb N} }

\def \R { {\mathbb R} }
\def \C { {\mathbb C} }
\def \Z { {\mathbb Z} }
\def \P { {\mathbb P} }

\def \S {{\mathcal S}^{*}(a,b)}
\def \Zz {{\mathcal Z}}

\DeclareMathOperator{\prob}{\P rob}

\DeclarePairedDelimiter{\floor}{\lfloor}{\rfloor}

\newcommand{\new}[1]{#1}
 \title{Exact Goodness-of-Fit Testing for the Ising Model}

\author[Mart\'in del Campo]{Abraham Mart\'in del Campo}
\address{CONACYT Research Fellow -- CIMAT\\
       Jalisco S/N, Valenciana\\
       36023 Guanajuato, Gto. \\
       M\'exico}
\email{abraham.mc@cimat.mx}
\urladdr{http://www.cimat.mx/~abraham.mc/}

\author[Cepeda]{Sarah Cepeda}
\address{IST Austria\\
         Am Campus 1\\
         3400 Klosterneuburg\\
         Austria}
\email{sarah.cepeda@ist.ac.at}

\author[Uhler]{Caroline Uhler}
\address{MIT, 77 Massachusetts Ave., 32-D634, Cambridge, MA 02139, USA \newline IST Austria\\
         Am Campus 1\\
         3400 Klosterneuburg\\
         Austria}
\email{cuhler@mit.edu}
\urladdr{http://www.carolineuhler.com/}

\date{\today}
\begin{document}

\begin{abstract}
The Ising model is one of the simplest and most famous models of interacting systems. It was originally proposed to model ferromagnetic interactions in statistical physics and is now widely used to model spatial processes in many areas such as ecology, sociology, and genetics, usually without testing its 
\new{goodness of fit}. 
Here, we propose \new{various test statistics and} an exact goodness-of-fit test for the finite-lattice Ising model. The theory of Markov bases has been developed in algebraic statistics for exact goodness-of-fit testing using a Monte Carlo approach. However, 
finding a Markov basis is often computationally intractable.
Thus, we develop a Monte Carlo method for exact goodness-of-fit testing for the Ising model which avoids computing a Markov basis and also leads to a better connectivity of the Markov chain and hence to a faster convergence. We show how this method can be applied to analyze the spatial organization of receptors on the cell membrane.
\end{abstract}

\maketitle

\section{Introduction}

The Ising model was invented as a mathematical model of ferromagnetism in statistical mechanics and it first appeared in~\cite{Ising}, a paper based on Ising's Ph.D.~thesis. The model consists of binary variables, called \emph{spins}, which are  usually arranged on an integer lattice, allowing each spin to interact only with its neighbors. The Ising model has played a central role in statistical mechanics; see \cite{Ising_history} for a review. Since the Ising model allows a simplified representation of complex interactions, it has 
been used in various areas, e.g.~for image analysis pioneered by Besag~\cite{Besag_image1, Besag_image2}, to model social and political behavior~\cite{Galam_1, Galam_2}, to model the interactions of neurons~\cite{Schneidman2006, TMASBB14, Tkacik}, or to analyze genetic data~\cite{Majewski2001}. Statistical inference for the Ising model is often made on the basis of a single observed lattice configuration. In this paper, we propose a method for goodness-of-fit testing for the Ising model that can be applied in this setting.

Among the many contributions of Julian Besag, he wrote a series of papers on hypothesis testing for spatial data when asymptotic approximations are inadequate~\cite{ Bes74, Bes77, BC89, BB00}. He proposed and applied Monte Carlo tests for this purpose~\cite{Bes72, Bes74}: An irreducible, aperiodic Markov chain starting in the observed configuration is built on the set of all spatial configurations with the same sufficient statistics. The true (non-asymptotic) conditional $p$-value is approximated by the $p$-value resulting from the distribution of a test statistic evaluated at the Monte Carlo samples. The main difficulty of this approach is to guarantee irreducibility of the Markov chain. In~\cite{BC89} Besag and Clifford discuss hypothesis testing specifically for the Ising model and propose using simple swaps of two randomly chosen lattice points of different states, accepting those which preserve  the sufficient statistics. As noted by Bunea and Besag in~\cite{BB00}, this algorithm leads to a \emph{reducible} Markov chain without any guarantees of converging to the correct conditional distribution. It remained an open problem to develop an irreducible version of this algorithm, a problem we solve in this paper.

Diaconis and Sturmfels introduced \emph{Markov bases} and developed a general framework for sampling from a conditional distribution (given the sufficient statistics) in discrete exponential families~\cite{DS98}. A Markov basis is a set of moves that connect any two configurations with the same sufficient statistics by passing only through configurations that preserve the sufficient statistics. Thus, a Markov basis allows building an irreducible Markov chain and performing non-asymptotic goodness-of-fit testing. As shown by Diaconis and Sturmfels in~\cite{DS98}, finding a Markov basis for a particular model is equivalent to finding generators for a specific ideal in a polynomial ring and can (in principle) be computed using Gr\"obner bases techniques. This spurred a lot of research in algebraic statistics (see~\cite{AHT12,DSS09,HoSu07,Rap03,PP13} or more recently~\cite{cai2014,GPS14,KOT15,RS16,SZP14}).
\new{However, computing a Markov basis is often computationally intractable, making the Markov basis approach to goodness-of-fit testing 
ineffective in many applications.}
As we show in Section~\ref{ss:MB}, a Markov basis for the Ising model on a $3{\times} 3$ lattice consists already of 1,334 moves. Computing a Markov basis for the Ising model on a $4{\times} 4$ lattice is computationally infeasible using the current Gr\"obner basis technology. A similar problematic has been observed 
\new{for other models, for instance for network models~\cite{PRF10}. Here the largest models considered so far are food webs of around 40 nodes~\cite{GPS14, OHT13}.} 
See~\cite{MBD} for a collection of Markov bases that have been computed during the last years.

In this paper, we propose a method for sampling from a conditional distribution (given the sufficient statistics) which circumvents the need of computing a Markov basis. The idea is to use a set of simple moves, a subset of the Markov basis, which can be computed easily but does not necessarily lead to a connected Markov chain. We then build a connected Markov chain using these simple moves by allowing \new{the sufficient statistics to change slightly. For fast convergence of the Markov chain it is desirable to determine the minimal amount of change needed in the sufficient statistics to guarantee connectedness and we prove that just allowing the direct neighbors is sufficient. This results in an irreducible version of the algorithm by Bunea and Besag~\cite{BB00} and we show how it can be applied for goodness-of-fit testing in the Ising model.}

\new{We briefly review other methods that have previously been proposed in oder to bypass the computational bottleneck of computing the whole Markov basis for goodness-of-fit testing and compare these methods to our approach: For applications to contingency tables it has been shown that a simpler subset of the full Markov basis is sufficient to build a connected Markov chain in certain restricted settings, with constraints either on the sufficient statistics~\cite{CDY10} or on the cell entries~\cite{RY10}. Note that this is possible since using a Markov basis to build a Markov chain guarantees its connectedness for \emph{any} sufficient statistic and one can do better for \emph{particular} sufficient statistics. A different approach is to use a \emph{lattice basis} instead of the full Markov basis~\cite[Chapter 1]{DSS09}. A lattice basis is a subset of a Markov basis and therefore easier to compute. However, in order to guarantee connectedness of the Markov chain, one needs to use integer combinations of the elements of a lattice basis~\cite[Chapter 16]{AHT12}. A third approach is based on computing the Markov basis dynamically; at each step, this framework finds a subset of the Markov basis elements that connect the current iterate to all its neighbors~\cite{Dob12, GPS14}. The advantage of our method compared to the first approach is that it is applicable to \emph{any} sufficient statistic. In comparison to the second and third approach, our method reduces the size of the moves even further and allows application to lattices of size $800 \times 800$ as required for the biological data set in Section~\ref{s:biodata}.
}




Our paper is structured as follows: In Section~\ref{s:background}, we introduce notation for the finite-lattice Ising model and describe the sufficient statistics for this model. We also introduce Markov bases in the context of Ising models and determine 
a Markov basis for the 2-dimensional Ising model of size $3{\times} 3$. Section~\ref{s:MkvChains} contains the main results of this paper. We prove that simple swaps are sufficient for constructing an irreducible, aperiodic, and reversible Markov chain when allowing a bounded change in the sufficient statistics of size $2^{d-1}$, where $d$ denotes the ambient dimension of the Ising model. In Section~\ref{s:TestIsing} 
we discuss several test statistics for the Ising model and we analyze their performance using simulated data in Section~\ref{s:simulations}. Finally, in Section~\ref{s:biodata} we apply our exact goodness-of-fit test for the Ising model to biological data and analyze the spatial organization of receptors on the cell membrane.

\section{Background}\label{s:background}

\subsection{Ising model}
The Ising model was originally introduced to study magnetic phase transitions and is now one of the most famous models of interacting systems~\cite{Ising_history, Ising}. The Ising model consists of a collection of binary random variables, the spins, which are usually arranged on an integer lattice, with the edges representing interaction between
 spins. In this section, we introduce some notation and recall the mathematical definition of the Ising model.

Let $[k]$ denote the set of integers $\{1,\dotsc, k\}$.
For $N_1,\dotsc, N_d\in \Z_{>0}$, let $L$ denote the \emph{$d$-dimensional lattice graph} of size $N_1\times \cdots \times N_d$, whose set of vertices $V$ consists of the elements in 
$[N_1]\times \dotsb \times [N_d]$
and whose edge set $E$ consists of pairs $i,j\in V$ such that $i$ and $j$ agree in all coordinates except for one in which they differ by one. To each vertex $i\in V$ we associate a binary random variable $Y_i$ taking values in $\{0,1\}$ and edges represent interaction between the adjacent random variables.
A \emph{configuration} is an element $y\in S:=\{ 0,1\}^{N_1 \times\dotsb \times N_d}$. To each configuration $y$, we associate the following two quantities:
\begin{equation}\label{eq:suff_stats}
T_1(y) := \sum_{i\in V} y_i \, , \qquad T_2(y) := \sum_{(i,j)\in E} |y_i- y_j|,
\end{equation}
where the first quantity counts the number of ones in the configuration, and the second counts the number of edges between vertices with different values.

\begin{figure}[tb]
\begin{center}
	\includegraphics[width=.17\textwidth]{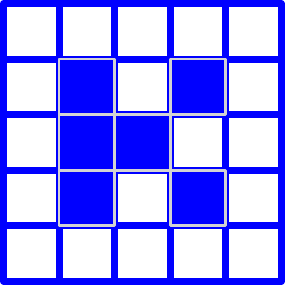}
\end{center}
\caption{A lattice configuration for the 2-dimensional Ising model.\label{F:grid1}}
\end{figure}

For $d=2$, we represent a configuration $y$ by a diagram, where the vertices of $L$ are represented by squares in a grid, the edges correspond to two 
adjacent squares, the ones in the configuration are represented by colored squares and the zeroes by white squares. For instance, the configuration 
\[
y = (
0,0,0,0,0, \,
0,{\color{blue}1},{\color{blue}1},{\color{blue}1},0, \,
0,0,{\color{blue}1},0,0, \,
0,{\color{blue}1},0,{\color{blue}1},0, \,
0,0,0,0,0
) \in \{0,1\}^{5\times 5}
\]
(the entries are read column-wise) is represented by the diagram depicted in Figure~\ref{F:grid1}. This example has $T_1(y) = 6$ and $T_2(y)=18$.

For the Ising model, the probability of observing a configuration $y\in S$ is given by the \emph{Boltzmann distribution}:
\begin{equation}\label{eq:Boltzmann}
P(y) = \frac{\exp(\alpha\cdot T_1(y)+\beta\cdot T_2(y))}{Z(\alpha,\beta)},
\end{equation}
where $\alpha, \, \beta$ are the model parameters and $Z(\alpha, \beta)$ is the normalizing constant. Note that the quantities $T_1(y)$ and $T_2(y)$ are the sufficient statistics of the model. Sometimes one considers an extra parameter representing the temperature of the configuration. In this paper, however, we focus on the simplified version where the joint probability does not depend on the temperature. 

Working over the state space $\tilde{S}:=\{ -1,1\}^{N_1 \times\dotsb \times N_d}$ and replacing every zero in a configuration $y\in S$ by $-1$, leads to a configuration $x\in\tilde{S}$, and we get the classical definition of the Boltzmann distribution with sufficient statistics 
\[
\widetilde{T}_1(x) := \sum_{i\in V} x_i, \qquad \widetilde{T}_2(x) := \sum_{(i,j)\in E} x_ix_j.
\]
Since
\begin{equation}\label{eq_relation}
\widetilde{T}_1(x) = 2\cdot T_1(y)-|V|,  \qquad \widetilde{T}_2(x) = |E| - 2\cdot T_2(y),
\end{equation}
where $|V|=\prod_{i=1}^d N_i$ is the number of vertices and $|E| = \sum_{i=1}^d ((N_i{-}1)\cdot \prod_{j\neq i} N_j)$ the number of edges in $L$, there is a one-to-one correspondence between the two representations. In the following we work over the state space $S$, since this simplifies notation.

Our goal is to develop an exact goodness-of-fit test for the Ising model that is applicable when only a single lattice configuration is observed. Our null hypothesis is that an observed configuration $y\in S$ is sampled from the Boltzman distribution~\eqref{eq:Boltzmann} for unspecified parameters $\alpha$ and $\beta$. In Section~\ref{s:TestIsing} we describe several alternatives such as the presence of long-range interactions or non-homogeneity.

For fixed $a,b \in \Z_{> 0}$ we define the corresponding \emph{configuration sample space} by
$$S(a,b)\;:=\,\{y\in S \;\mid\; T_1(y)=a,\; T_2(y)=b\}.$$
We denote the distribution of a Boltzmann random variable $Y$ conditioned on the values of the sufficient statistics by
\begin{equation}\label{eq:p-value}
\pi(y)\,:=\,\prob(Y=y \mid T_1(Y)=a,\; T_2(Y)=b).
\end{equation}
A simple computation shows that for the Ising model, $\pi$ is the uniform distribution on $S(a,b)$. Then the conditional $p$-value of a given configuration $y$ can be approximated using Markov chain Monte Carlo (MCMC) algorithms for sampling configurations from the conditional distribution~$\pi$. However, MCMC methods for hypothesis testing in the Ising model have various computational obstructions: the sample space $S(a,b)$ grows exponentially with the size of the lattice $L$ rendering direct sampling infeasible. In addition, to ensure that the conditional distribution $\pi$ is the stationary distribution of the Markov chain, the chain needs to be \emph{irreducible} (i.e.~connected). This is difficult to achieve since the 
space $S(a,b)$ has a complicated combinatorial structure. One way to overcome this problem is to use a Markov basis to construct the Markov chain. We now formally introduce Markov bases and then discuss the computational limitations of this approach 
in the Ising model. 

\subsection{Markov bases}\label{ss:MB}
First note that every integer vector $z\in \Z^{N_1\times \cdots \times N_d}$ can be written uniquely as a difference $z=z^+ - z^-$ of two nonnegative vectors with disjoint support. A Markov basis in the context of the Ising model is defined as follows.
\begin{defn}
A \emph{Markov basis} for the Ising model is a set $\widetilde{\mathcal Z}\subset \Z^{N_1\times \cdots \times N_d}$ of integer vectors such that:
\begin{enumerate}
\item[(i)] All $z\in \widetilde{\mathcal Z}$ satisfy 
$T_1(z^+) = T_1(z^-)$ and $T_2(z^+) = T_2(z^-)$.
\item[(ii)] For any $a,b \in \Z_{> 0}$ and any $x,y \in S(a,b)$, there exist $z_1,\ldots, z_k\in \widetilde{\mathcal Z}$ such that 
$$
  y = x+\sum_{i=1}^k z_i \quad \text{and } \quad x+\sum_{i=1}^\ell z_i \in S(a,b) \quad \mbox{for all } \ell = 1,\ldots, k.
  $$
\end{enumerate}
The elements of $\widetilde{\mathcal Z}$ are called \emph{moves}. A Markov basis allows constructing an aperiodic, reversible, and irreducible Markov chain that has stationary distribution $\pi$ using the Metropolis-Hastings algorithm, by selecting a proposed move $z$ uniformly from 
$\widetilde {\mathcal Z}$
and accepting the move from $y$ to $y+z$ with probability $\min(1, \pi(y+z)/\pi(y))$ as long as $y+z\in S(a,b)$. See~\cite[Lemma 2.1]{DS98} for a proof. Since 
$\pi$ is the uniform distribution for the Ising model, the acceptance probability is 1 for any element $y+z\in S(a,b)$.
\end{defn}

As explained in~\cite[Section 3]{DS98}, a Markov basis can also be defined algebraically. For the Ising model the algebraic description is as follows: consider the polynomial ring $\C[P_y \mid y\in S]$ of complex coefficients where the indeterminates are indexed by configurations. Similarly, we consider the polynomial ring $\C[u_i, v_i, p_{ij}, q_{ij} \mid i\in V,\, {(i,j)}\in E ]$ with two sets of indeterminates, one indexed by the vertices of $L$ and the other by the edges of $L$.
Given a configuration $y\in S$, let $\mathcal{L}_y = (\mathcal{V}_y, \mathcal{E}_y)$ denote the induced subgraph of $L$ defined by restricting $L$ to the vertex set $\mathcal{V}_y:= \{ i\in V \mid y_i =1\}$. In addition, let $E_y := \{ (i,j)\in E \mid y_i \neq y_j\}$.
Then we define a monomial map
\begin{equation}\label{eq:toricmap1}
\phi : \C[P_y] \longrightarrow \C[u_i, v_i, p_{ij}, q_{ij}], \qquad \phi(P_y) :=
\prod_{i\in \mathcal V_{y}} v_i \prod_{i\in V\setminus\mathcal V_{y}} u_i 
\prod_{(i,j)\in E_y} p_{ij} \prod_{(i,j)\in E \setminus E_{y}} q_{ij}.
\end{equation}
Let  $\C[u,v,p,q]$ be the polynomial ring in four indeterminates and consider the monomial ring homomorphism
\begin{equation}\label{eq:toricmap2}
\psi : \C[u_i, v_i, p_{ij}, q_{ij}] \longrightarrow \C[u,v,p,q]
\end{equation}
defined by 
$$
\psi(u_i) := u,\quad \psi(v_i) := v, \quad \psi(p_{ij}):= p,\quad \psi(q_{ij}):=q, 
$$
for all $i\in V$ and $(i,j)\in E$. With these definitions we obtain that 
$$
(\psi\circ \phi)(P_y) = v^{T_1(y)} u^{|V|-T_1(y)} p^{T_2(y)} q^{|E|-T_2(y)},
$$
where, as in (\ref{eq_relation}), $|V|=\prod_{i=1}^d N_i$ is the number of vertices and $|E| = \sum_{i=1}^d ((N_i{-}1)\cdot \prod_{j\neq i} N_j)$ is the number of edges in $L$. Then it follows directly from~\cite[Section 3]{DS98} that any generating set for the ideal $I:= \ker (\psi\circ \phi)$ is a Markov basis for the Ising model.
\begin{ex}\label{Ex:Markov}
Let $y = (0,1,0,0,1,1,0,0,0)$ be a configuration on a $3\times 3$ lattice $L$ as represented by the diagram 
\[
\vcenteredinclude{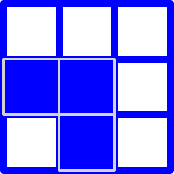}{35pt} \ .
\]
We denote the vertices of $L$ column-wise by $\{1,2,\ldots, 9\}$. Then the induced subgraph $\mathcal{L}_y$ has vertices $\mathcal V_y = \{2,5,6\}$ corresponding to the colored squares in the diagram. In this example, $E_y = \{(1,2), (2,3), (3,6),(4,5), (5,8), (6,9)\}$. Therefore, the image of $y$ under $\phi$ is
$$
\phi(P_y) = v_2\, v_5\, v_6\, u_1\, u_3\, u_4\, u_7\, u_8\, u_9\, p_{12}\, p_{23}\, p_{36}\, p_{45}\, p_{58}\, p_{69}\, q_{14}\ q_{25}\, q_{47}\, q_{56}\, q_{78}\, q_{89},
$$
and the image under the composition of $\phi$ and $\psi$ is
$$
(\psi\circ \phi)(P_y) = v^3 u^6 p^6 q^6,
$$
where the sufficient statistics $T_1(y)$ and $T_2(y)$ are the exponents of the indeterminates $v$ and $p$, respectively. \qed 
\end{ex}

For the $3{\times}3$ lattice, we computed a Markov basis for the ideal $I$ using the algebraic software {\tt 4ti2}~\cite{4ti2}. The Markov basis consists of 466 binomials of degree one, 864 binomials of degree two, and 4 binomials of degree three. For instance, the following degree one generator of $I$ involves the configuration $y$ from Example~\ref{Ex:Markov}:
$$
P_{\, \includegraphics[height=13pt]{Pictures/333_1.pdf}} - P_{\, \includegraphics[height=13pt]{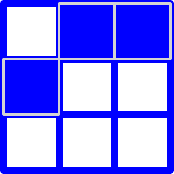}}.
$$
Note that both configurations have the same sufficient statistics. Sampling this move would mean that if the chain is in the configuration on the right, we would move to the configuration on the left.
Examples of elements of degree two and three in the computed Markov basis are:
$$
P_{\, \includegraphics[height=13pt]{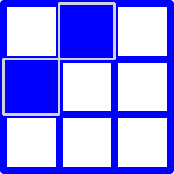}} 
P_{\, \includegraphics[height=13pt]{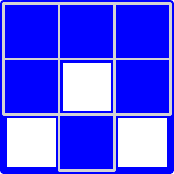}}
- P_{\, \includegraphics[height=13pt]{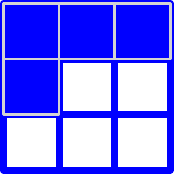}} 
P_{\, \includegraphics[height=13pt]{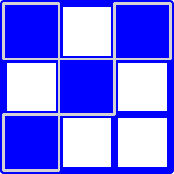}}
\qquad \mbox{and }\quad 
P_{\, \includegraphics[height=13pt]{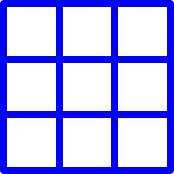}}
P_{\, \includegraphics[height=13pt]{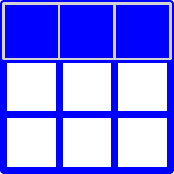}}
P_{\, \includegraphics[height=13pt]{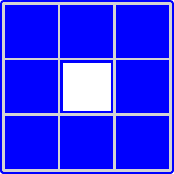}}
- P_{\, \includegraphics[height=13pt]{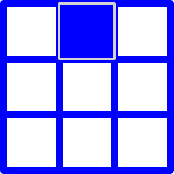}}
P_{\, \includegraphics[height=13pt]{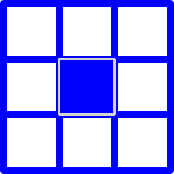}}
P_{\, \includegraphics[height=13pt]{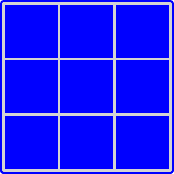}}\ .
$$
Note that the sum of the sufficient statistics of the configurations on the positive and the negative side are the same. When hypothesis testing is based on one observed configuration only, then only the degree one moves are of interest; however, if we had multiple observed configurations, then we would also use the higher degree moves. Since the number of indeterminates involved in these Gr\"obner bases computations grows exponentially with the size of the lattice, we were not able to compute a Markov basis for a lattice of size $4{\times} 4$ or larger using the current technology. One could optimize the Markov basis computation by considering 
\new{other parameterizations of the model (see~\cite{PR12}) or by computing a Markov subbasis~\cite[Chapter 16]{AHT12} that contains only the degree 1 moves.  Also, it would be interesting to study the structure of the Universal Gr\"obner basis or the Graver basis~\cite[Chapter 1]{DSS09} (both contain the Markov basis as a subset) for the Ising model, since these bases are often easier to describe algebraically, in particular in the presence of symmetries. However, even with these approaches it would be difficult to compute the moves needed for a lattice of size $10\times 10$ (requires computations in a ring with over 100 indeterminates), let alone a lattice of size $800\times 800$, as would be required for the example discussed in Section~\ref{s:biodata}.}


\section{Constructing an irreducible Markov chain}\label{s:MkvChains}

In~\cite{BC89}, Besag and Clifford describe a Metropolis-Hastings algorithm for conditional hypothesis testing in the Ising model. At each iteration two lattice points of different states are randomly selected and their states are switched if this preserves the sufficient statistics. However, as noted by Bunea and Besag in~\cite{BB00}, this algorithm in general produces a reducible Markov chain on $S(a,b)$; hence, it may not result in uniform samples from $S(a,b)$ and lead to inaccurate $p$-values. On the other hand, as we have seen in the previous section, computing a Markov basis for this problem guarantees irreducibility but is computationally intractable. In the following, we show how to modify the algorithm in~\cite{BC89} to obtain an irreducible Markov chain without computing a Markov basis. Our method consists of using the simple moves proposed by Besag and Clifford, but to expand the sample space and allow 
\new{the sufficient statistics to change slightly.} 
In this section we prove that the resulting Markov chain is irreducible by exploiting the combinatorics underlying the Ising model.

We define a \emph{simple swap} to be an integer vector $z\in \Z^{N_1\times \cdots \times N_d}$ of the form $z = {\bf e}_i - {\bf e}_j$, where the ${\bf e}_i$'s denote the canonical basis vectors of $\Z^{N_1\times \cdots \times N_d}$. Simple swaps correspond to switching the states of two lattice points in configuration $y$, i.e., replacing a pair $(y_i{=}0, y_j{=}1)$ by $(y_i{=}1, y_j{=}0)$. Throughout the paper, $\mathcal Z$ denotes \emph{the set of simple swaps}. Let $y\in S(a,b)$ be a configuration with $y_i {=} 0$ and $y_j{=}1$ and let $z\in\mathcal{Z}$ be the simple swap ${\bf e}_i - {\bf e}_j$. Then $T_1(y) = T_1(y{+}z)$ by construction, but $T_2(y+z)$ may disagree with $T_2(y)=b$. We say that two configurations $y, y' \in S(a,b)$ are $S(a,b)$-\emph{connected} by $\mathcal Z$, if there is a path between $y$ and $y'$ in $S(a,b)$ consisting of simple swaps $z\in \mathcal Z$, i.e., if there exist $z_1,\dots , z_k\in\mathcal{Z}$ such that
$$y'=y+\sum_{i=1}^k z_i, \quad \textrm{with}\quad y+\sum_{i=1}^{\ell} z_i\in S(a,b) \quad\textrm{ for all }\, \ell=1,\dots ,k.$$ 

With this notation the algorithm proposed by Besag and Clifford~\cite{BC89} is as follows: start the Markov chain in a configuration $y\in S(a,b)$. Select $z\in \mathcal Z$ uniformly at random and let $y'=y{+}z$. Accept $y'$ if $y' \in S(a,b)$; otherwise, remain in $y$. \new{We illustrate with the following example that such a chain might not be able to leave its initial state.}
%

\begin{ex}\label{ex:NotConnected}
Consider the two configurations $y, y'$ in $\{0,1\}^{4\times 6}$ shown in Figure~\ref{F:counterexample}. Note that $y, y' \in S(4,8)$ but 
any simple swap $z \in \mathcal{Z}$ satisfies $T_2(y {+} z)> T_2(y)$. Therefore, 
the two configurations $y,y'$ are not $S(4,8)$-connected by $\mathcal{Z}$,
although they have the same sufficient statistics.\qed
\begin{figure}[htb]
\begin{tabular}{rcl}
	\vcenteredinclude{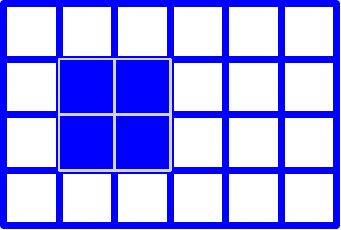}{.09\linewidth}
	& $\longleftrightarrow$ & 
	\vcenteredinclude{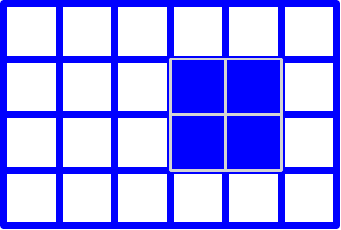}{.09\linewidth}
\end{tabular}
\caption{Simple swaps may not suffice to connect two configurations with the same sufficient statistics.\label{F:counterexample}}
\end{figure}
\end{ex}

Interestingly, as we show in Section~\ref{s:1dIsing}, this problem does not arise for the 1-dimensional Ising model. In higher dimensions we overcome this problem by running the Metropolis-Hastings algorithm in the 
smallest \emph{expanded sample space} $\S$, a set
with  $S(a,b)\subseteq \S$ such that any two configurations $y,y' \in S(a,b)$ are $\S$-connected by $\mathcal Z$, and then project the resulting Markov chain to $S(a,b)$. The resulting Markov chain is also reversible and aperiodic, since there is some holding probability. Before presenting the construction of the set $\S$ and the proof of irreducibility of the resulting Markov chain, we introduce a special class of configurations.

\subsection{Max-singleton configurations} 

A crucial element in the proof of irreducibility consists of configurations that maximize the number of singletons. Let $L=(V,E)$ be the graph of a $d$-dimensional integer lattice of size $N_1{\times}\cdots \times N_d$. Then each vertex $i\in V$ is adjacent to $2 d$ other vertices, except for those on the boundary of $L$. To simplify our exposition and not worry about boundary effects, \emph{we assume that all configurations have value 0 on the boundary}.

Let $y\in S$ be a configuration and $\mathcal{L}_y$ the induced subgraph defined by the vertices $\mathcal{V}_y$. For a connected component $C$ of $\mathcal L_y$, the \emph{size} $|C|$ is the number of vertices in $C$. A connected component of size one is called a \emph{singleton}. A \emph{max-singleton configuration} for $S(a,b)$ is a configuration $y^*\in S(a,b)$ that maximizes the number of singletons.

\subsection{Irreducibility in the 1-dimensional Ising model}\label{s:1dIsing}

We now prove that the original algorithm by Besag and Clifford~\cite{BC89} using simple swaps results in an irreducible Markov chain for the 1-dimensional Ising model. Recall that the set $\S$ is defined as the smallest set that contains the sample space $S(a,b)$ and for which any two configurations in $S(a,b)$ can be connected by simple swaps $\Zz$ without leaving $\S$. So we need to show that $\S=S(a,b)$. We do this by proving that the max-singleton configuration in $S(a,b)$ for the 1-dimensional Ising model is unique and then showing that any configuration $y\in S(a,b)$ is $S(a,b)$-connected by simple swaps to the unique max-singleton configuration. The proofs of the following results are given in the appendix.

\begin{lem}\label{L:shapeOptimal1D}
The max-singleton configuration of $S(a,b)$ for the 1-dimensional Ising model is unique (up to location
of its connected components) and consists of $b/2-1$ singletons and one connected component of size $a-b/2+1$.
\end{lem}
\begin{thm}\label{P:irred1D}
For $a, b\in \N$ and $y\in S(a,b)$, let $y^*\in S(a,b)$ denote the unique max-singleton configuration. There exists a sequence $z_1,\ldots, z_k\in \Zz$ such that 
\[
 y^* = y+\sum_{i=1}^k z_i \quad \text{and } \quad  y+\sum_{i=1}^\ell z_i \in S(a,b) \quad \mbox{for } \ell = 1,\ldots, k.
\]
\end{thm}

By noting that all steps in the proof of \new{Theorem}~\ref{P:irred1D} are reversible, we obtain the irreducibility of the Markov chain based on simple swaps for the 1-dimensional Ising model.
%
%


\subsection{Irreducibility in the 2-dimensional Ising model}\label{s:2dIsing}

Example~\ref{ex:NotConnected} illustrated that simple swaps are not sufficient for constructing an irreducible Markov chain for the 2-dimensional Ising model. Notice that for any configuration $y\in S(a,b)$, the number $b=T_2(y)$ is always even (assuming that $y$ has zeros on the boundary). This can be seen by analyzing the induced graph $\mathcal{L}_y=(\mathcal{V}_y, \mathcal{E}_y)$ and noting that $b=4a-2|\mathcal{E}_y|$. Hence, since simple swaps preserve the sufficient statistic $T_1(y)$, the minimal expansion of $S(a,b)$ is given by $S(a,b)\cup S(a,b\pm2)$. We prove in this section that this minimal expansion is sufficient for the 2-dimensional Ising model, i.e. $\S = S(a,b)\cup S(a,b\pm2)$. As for the 1-dimensional case, we start by identifying the max-singleton configuration and then show that any configuration $y\in S(a,b)$ is $\S$-connected by $\Zz$ to the max-singleton configuration. The technical details can be found in the appendix.

\begin{defn}\label{D:rectangularConfig}
A \emph{rectangular configuration} $y_\square$ is defined by $(n,m,d_1,d_2,s)$ with $n,m,d_1, d_2, s\in \Z_{>0}$, $d_2\in\{0,1\}$, $d_2\leq d_1< m\leq n$, and is of the following form: The induced subgraph $\mathcal{L}_y$ consists of $s$ singletons and one additional connected component of size $nm{+}d_1{+}d_2$, namely a rectangle $B$ of size $n{\times} m$ with a block of size $d_1{\times} 1$  lying on one of the shorter sides of $B$ and $d_2$ vertices on the longer side of $B$.
\end{defn}

\begin{figure}[htb]
	\includegraphics[width=.30\textwidth]{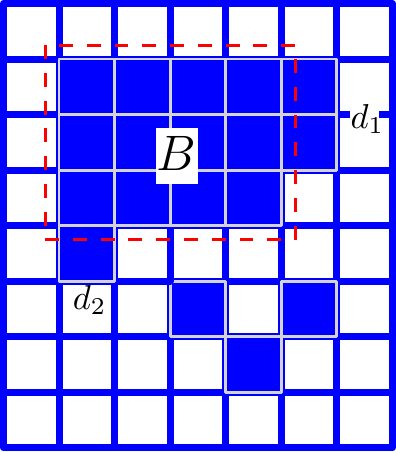}
\vspace{-5pt}
\caption{A rectangular configuration $y_{\square}=(3,4,2,1,3)$.\label{F:rectangular}}
\end{figure}
Figure~\ref{F:rectangular} 
shows a rectangular configuration, where $B$ is a rectangle of size $3{\times} 4$, i.e. $n=4$, $m=3$, and $d_1=2$, $d_2=1$, and $s=3$. Note that a rectangular configuration $y_\square$ defined by $(n,m,d_1, d_2, s)$ has sufficient statistics 
\begin{equation}\label{Eq:SuffStats}
T_1(y_\square) = nm + d_1+d_2 +s, \qquad 
T_2(y_\square) = 2(n+m+\delta_{d_1}+d_2+2s),
\end{equation}
where $\delta_{d_1}=0$ when $d_1=0$, and $\delta_{d_1} =1$ otherwise. The following result shows that the max-singleton configuration for the 2-dimensional Ising model is a rectangular configuration. The proof is given in the appendix.

\begin{lem}\label{lemma_opt_2d}
The max-singleton configuration of $S(a,b)$ for the 2-dimensional Ising model is unique (up to location) and corresponds to a rectangular configuration $(n,m,d_1,d_2,s)$ with 
\begin{align*}
&s = \Bigl\lfloor\dfrac{b/2-1-\sqrt{4a-b+1}}{2}\Bigr\rfloor, \quad 
m = \floor{\sqrt{a-s}},\quad 
r = \Bigl\lfloor\dfrac{a-s}{m}\Bigr\rfloor,\\ 
& n = \begin{cases}
r-1      & \text{if } \frac{a-s}{m}=\lfloor \frac{a-s}{m}\rfloor \text{ and } b-4s -2(m+r)=2, \\
  r    & \text{otherwise},
\end{cases}\\
&(d_1, d_2) =  \left\{ \begin{array}{ll} (a-s-mn, 0) & \textrm{if } b-4s-2(m+r)= 0\\ (a-s-mn-1, 1) & \textrm{if } b-4s-2(m+r) =2\end{array}. \right.
\end{align*}
\end{lem}

As an example consider the configurations of $S(18,30)$.  The max-singleton configuration prescribed by Lemma~\ref{lemma_opt_2d}
is a rectangular configuration with $s=3$, $m=3$, $n=4$, $d_1 =2$, and $d_2=1$. This is precisely the one shown in Figure~~\ref{F:rectangular}.

To prove that $\S = S(a,b)\cup S(a,b\pm2)$ we first show that any rectangular configuration $y_\square\in S(a,b)$ is $\S$-connected by $\Zz$ to the max-singleton configuration $y^*\in S(a,b)$. Then we show that any configuration $y\in S(a,b)$ is $\S$-connected by $\Zz$ to a rectangular configuration $y_\square\in S(a,b)$. By reversing the steps, this proves that any two configurations $y,y'\in S(a,b)$ are $\S$-connected by $\Zz$. The proof of the following result is simple but technical, and is therefore given in the appendix.

\begin{prop}\label{P:rectangle2square}
Let $y_\square\in S(a,b)$ be a rectangular configuration and let $x^*\in S(a,b)$ be the max-singleton configuration. Then there exits a sequence $z_1,\ldots, z_k\in \mathcal Z$ such that 
$$
x^* = y_\square+\sum_{i=1}^k z_i \quad \mbox{and }\quad y_\square+\sum_{i=1}^\ell z_i \in \S \quad \mbox{for all }\,\ell =1,\ldots, k.
$$
\end{prop}

We conclude with our main connectivity result which proves irreducibility of the Metropolis-Hastings algorithm. To simplify the argument involved in the proof, we assume that we have enough space to locate configurations across the lattice $L=(V,E)$, i.e., that $a << |V|$.

\begin{thm}
Assume $a<< |V|$. For any $y,y'\in S(a,b)$, there exist $z_1,\ldots, z_k\in \mathcal Z$ such that 
$$
y = y'+\sum_{i=1}^k z_i \quad \mbox{and }\quad y'+\sum_{i=1}^\ell z_i \in \S \quad \mbox{for all }\,\ell =1,\ldots, k.
$$
\end{thm}
\begin{figure}[tb]
\centering
\subfigure[Bounding region and its outermost layer.]
{\raisebox{3pt}{
\includegraphics[scale=1.06]{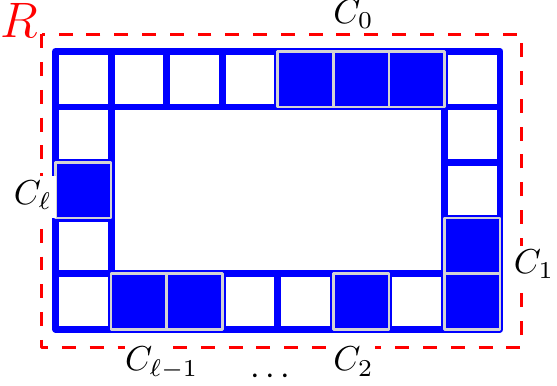}
}}
\label{F:bounding_region}
\subfigure[The rectangular configuration $\widetilde C$.]
{
\includegraphics[scale=1.3]{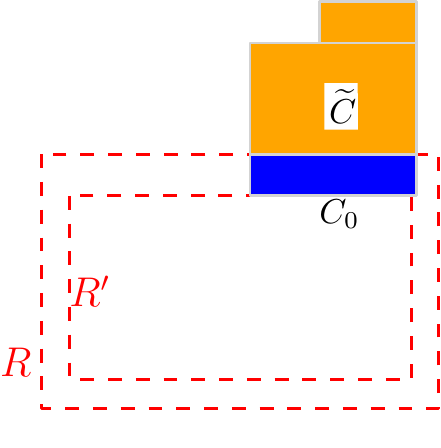}
}
\label{F:rectangular_Ctilde}
\caption{Construction of the rectangular configuration $\widetilde C$.}
\label{F:boundingR}
\end{figure}

\begin{proof}
As a consequence of Proposition~\ref{P:rectangle2square}, it suffices to prove that any configuration $y\in S(a,b)$ is $S^*(a,b)$-connected to a rectangular configuration $y_\square\!\!\in\!S(a,b)$. Let $R$ denote the smallest rectangle containing $\mathcal{L}_y$. To simplify the argument, we assume that $L>>R$. First, we remove all singletons from $R$ and place them in $L\setminus R$. Then we denote by $C_0, \ldots, C_\ell$ the connected components in the outermost layer of $R$ as illustrated in Figure~\ref{F:boundingR} (a) and we assume without loss of generality that $C_0$ is the largest such connected component. 

We now show that by simple swaps within $S^*(a,b)$, we can move the connected components $C_1,\dots, C_\ell$ to form a rectangle $\widetilde{C}$ with basis $C_0$ lying outside of $R$, as depicted in Figure~\ref{F:boundingR}~(b).  We move each connected component $C_i$ by repeatedly swapping the extremal sites of smallest degree. 
In Figure~\ref{F:possible moves}  we
illustrate the different scenarios and the corresponding change in $T_2(y)$ for moving a vertex (indicated by a dashed box) in table A to the position indicated in table B. 
If $|C_i| = 1$, then it is either a singleton and we move it outside of $R$, or it is a vertex of degree 1 (see A1). Starting a new row on top of $C_0$ corresponds to scenario B2. So the first move is a 0-swap (A1$\to$B2), meaning that $T_2(y)$ remains the same. If $|C_i|>1$, then the extremal sites are of degree 1 (see A2) or of degree 2 (see A3). So in this case the first move is either a 0-swap (A2$\to$B2) or a $+2$-swap (A3$\to$B2). The remaining swaps to move $C_i$ on top of $C_0$ are either $-2$-swaps (A1$\to$B1 or A2$\to$B1) or $0$-swaps (A3$\to$B1). If the configuration lies in $S(a,b{-}2)$ and the next move is of type A1$\to$B1, we instead create a singleton (A1$\to$B3) and get a configuration in $S(a,b)$. Note that if the first move is a $+2$-swap (A3$\to$B2), then the last remaining move in $C_i$ is a $-2$-swap (A1$\to$B1), because we first swap extremal sites of smallest degree. So since $|C_0|\geq |C_i|$, we are never in the situation in which we make two $+2$-swaps before a $-2$-swap.

\begin{figure}[tb]
\begin{tabular}{|c|c|c|}
\hline
\multicolumn{3}{|c|}{A} \\ \hline
1          & 2      &3    \\ \hline
\raisebox{-2pt}{\includegraphics[width=.08\textwidth]{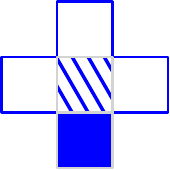} }       & 
\raisebox{-5pt}{\includegraphics[width=.08\textwidth]{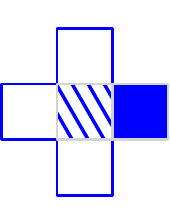}}        &
\raisebox{-2pt}{\includegraphics[width=.08\textwidth]{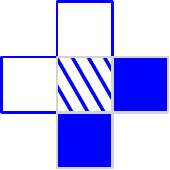} } 
\\[5pt] \hline
\end{tabular}
\qquad\qquad\begin{tabular}{|c|c|c|}
\hline
\multicolumn{3}{|c|}{B} \\ \hline
1          & 2     & 3     \\ \hline
\raisebox{-5pt}[25pt][0pt]{\includegraphics[width=.08\textwidth]{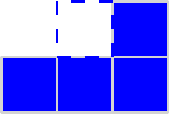} }   & 
\raisebox{-5pt}{\includegraphics[width=.08\textwidth]{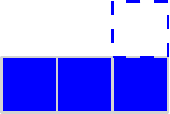}}    & 
\raisebox{0pt}{\hspace{5pt}\includegraphics[width=.026\textwidth]{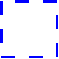}\hspace{5pt}}\\[5pt] \hline
\end{tabular}\bigskip
\\
\begin{tabular}{|c||c||c||c||c||c||c|}
\hline
 A1$\to$B1 & A1$\to$B2 & A1$\to$B3 & A2$\to$B1 & A2$\to$B2  & A3$\to$B1 & A3$\to$B2 \\\hline
$-2$-swap & $0$-swap & $+2$-swap & $-2$-swap & $0$-swap & $0$-swap & $+2$-swap  \\\hline
\end{tabular}
\vspace{10pt}
\caption{Types of simple swaps we use.}\label{F:possible moves}
\end{figure}

In this way, we transform the connected components $C_1,\dots, C_\ell$ into a rectangular configuration $\widetilde{C}$. We repeat the process by letting $C'_0,C_1' , \ldots, C'_{\ell'}$ be the components of the next layer in $R$ with $C'_0$ denoting the largest connected component. If $|C'_0| \leq |C_0|$, then we move $C'_0,C_1' , \ldots, C'_{\ell'}$ on top of $\widetilde{C}$ as before. Otherwise, we move $C_0\cup \widetilde{C}$ on top of $C'_0$ and proceed with $C_1' , \ldots, C'_{\ell'}$ as before with $C_1 , \ldots, C_{\ell}$. We continue inductively until we reach a rectangular configuration (with $d_2=0$) that lies either in $S(a,b)$ or $S(a,b{-}2)$. In the latter case, we obtain a rectangular configuration in $S(a,b)$ by reducing $d_1$ by one and letting $d_2=1$ as in the proof of Proposition~\ref{P:rectangle2square} in the appendix (see Figure~\ref{F:breakcorner}).
\end{proof}

\subsection{Higher dimensional Ising model}
We now show how our Metropolis-Hastings algorithm can be extended to the $d$-dimensional Ising model when $d\geq 3$. Let $L=(V,E)$ be the $d$-dimensional lattice graph of size $N_1\times \cdots \times N_d$. Let $\mathcal{S}^*_d(a,b)$ denote the smallest expanded sample space such that any $y,y'\in S(a,b)$ in the $d$-dimensional lattice are $\mathcal{S}^*_d(a,b)$-connected by $\Zz$. In Section~\ref{s:1dIsing} we showed that  $\mathcal{S}^*_1(a,b)=S(a,b)$ and in Section~\ref{s:2dIsing} we proved that $\mathcal{S}^*_2(a,b) = S(a,b)\cup S(a,b\pm 2)$. Since any $d{-}1$-dimensional configuration can be realized as a configuration in a $d$-dimensional lattice, it must hold that $\mathcal{S}^*_{d-1}(a,b)\subseteq \mathcal{S}^*_d(a,b)$. Moreover, let us consider the $d$-dimensional version of Example~\ref{ex:NotConnected}, where $y$ is a hyper-rectangular configuration of size $2\times 2\times \cdots \times 2$ in the $d$-dimensional lattice. Note that any simple swap $z\in\Zz$ satisfies $T_2(y{+}z)=T_2(y)+2d$ or  $T_2(y{+}z) = T_2(y)+2(d{-}1)$ and hence $S(a,b\pm 2(d{-1}))\subseteq \mathcal{S}^*(a,b)$. We prove that this expansion is in fact sufficient, i.e.
$$
\mathcal{S}^*_d(a,b) = \mathcal{S}^*_{d-1}(a,b) \cup S(a,b\pm 2(d-1)).
$$
For instance, for the 3-dimensional lattice this means that $\mathcal{S}^*_3(a,b) = S(a,b)\cup S(a,b\pm2) \cup S(a,b\pm 4)$. The proof of irreducibility in the $d$-dimensional case is analogous to the 2-dimensional case. We therefore omit the details and only show how to generalize the two most important elements of the proof, namely the concept of a \emph{rectangular configuration} and the types of simple swaps  in Figure~\ref{F:possible moves} to connect any configuration to a rectangular configuration. 

\begin{defn}\label{D:rectangularD}
A \emph{$d$-dimensional rectangular configuration} consists of $s$ singletons and one additional connected component $B$ that decomposes into $B= (B_d, B_{d-1}, \ldots, B_1, d_2)$, where $B_d$ is a $d$-dimensional rectangular block of size $m_1{\times}m_2\times \cdots \times m_d$ with $m_1\leq m_2\leq \cdots \leq m_d$. For each $k=2, \ldots , d{-}1$, $B_k$ is a $k$-dimensional rectangular configuration that lies in the smallest facet (of size $m_1\times \cdots \times m_{k-1}$) of the $k{+}1$-dimensional rectangular block $B_{k+1}$. Lastly, $d_2\in\{0,1\}$. 
\end{defn}

%
\begin{figure}[tb]
	\includegraphics[width=.25\textwidth]{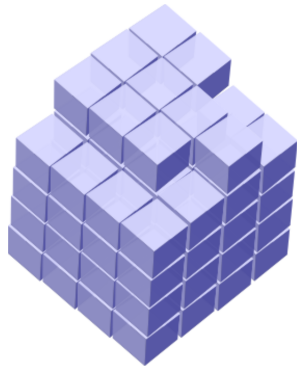}
\vspace{-5pt}
\caption{A 3-dimensional rectangular configuration.\label{F:rectangular3d}}
\end{figure}

A 3-dimensional rectangular configuration is shown in Figure~\ref{F:rectangular3d}. The $d$-dimensional analogue of the simple swaps in Figure~\ref{F:possible moves} are swaps of a vertex
that is adjacent to either $1, 2, \ldots, d$ neighbors (as indicated in Table A in Figure~\ref{F:possible moves} for the 2-dimensional case), to a position (indicated in Table B in Figure~\ref{F:possible moves} for the 2-dimensional case) with either  $0, 1, 2, \ldots, d{-}1$, or $d$ neighbors. With this observation it is straight-forward to generalize the proofs for the 2-dimensional setting to show irreducibility of our Metropolis-Hastings algorithm for the Ising model in higher dimensions.

\begin{thm}
Suppose $a<<|V|$. For any $d$-dimensional configurations $y,y'\in S(a,b)$, there exits a sequence $z_1,\ldots, z_k\in \mathcal Z$ such that 
$$
y = y'+\sum_{i=1}^k z_i \quad \mbox{and }\quad y'+\sum_{i=1}^\ell z_i \in \S \quad \mbox{for all }\, \ell =1,\ldots, k.
$$
\end{thm}


\section{Test statistics for the Ising model}\label{s:TestIsing}

After showing how to construct an irreducible, aperiodic and reversible Markov chain on $S(a,b)$, we now describe various statistics for testing for departure from the Ising model. The null model is the Ising model and departure from the null model is possible either by the presence of long-range interactions or by the presence of non-homogeneity. In the following, we describe both alternatives and various test statistics.

%
%
%
%

\subsection{Presence of long-range interactions}

This alternative hypothesis is defined by
\begin{equation}\label{E:alternative1}
\prob(y) = \frac{\exp(\alpha\cdot T_1(y)+ \beta\cdot T_2(y) + \gamma \cdot T_3(y))}{Z(\alpha, \beta, \gamma)},
\end{equation}
where $\alpha, \beta$, and $\gamma$ are the model parameters, $Z(\alpha, \beta, \gamma)$ is the normalizing constant, and $T_1(y)$, $T_2(y)$, and $T_3(y)$ are the minimal sufficient statistics, where $T_1(y)$ and $T_2(y)$ are the sufficient statistics of the Ising model as defined in~\eqref{eq:suff_stats}. For instance, $\gamma$ could correspond to the \emph{interaction between second-nearest neighbors} with
$$
T_3(y) \;=\; \sum_{(i,j), (j,k)\in E} |y_i-y_k|.
$$
A particular example of this model is a model with \emph{interaction between diagonal neighbors} with 
\begin{equation}\label{E:diagonal_interaction}
T_3(y) \;=\; \sum_{\substack{(i,j), (j,k)\in E \\ i,k \textrm{ diagonal neighbors}}} |y_i-y_k|.
\end{equation}
Or $\gamma$ could correspond to an \emph{overall effect} for example with 
\begin{equation}\label{E:overall_interaction}
T_3(y) \;=\; \left\{ \begin{array}{ll}
 1 & \textrm{if \;$\sum_{i\in V} y_i$\, is even,}\\
0 & \textrm{otherwise},
  \end{array} \right.
\end{equation}
an indicator function taking the value 1 if the number of ones in the lattice is even and taking the value 0 otherwise.

In ~\cite{Bes77}, Besag proposed as a test statistic for the Ising model to use
\begin{itemize}
\item the \emph{number of diagonal pairs}, i.e.~the number of patterns of the form \vcenteredinclude{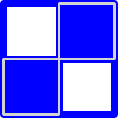}{15pt} (and its rotational analog) in a configuration.
\end{itemize}
Applied in a 1-sided test, this statistic may be an indicator for the presence of positive higher-order interaction. However, one also expects to see a large number of diagonal pairs under the Ising model with negative nearest-neighbor interaction. As we show in our simulations in Section~\ref{s:simulations}, other statistics are therefore usually more adequate for goodness-of-fit testing in the Ising model. We propose and analyze the following two test statistics to be used for 2-sided tests, namely
\begin{itemize}
\item the \emph{number of adjacent pairs}, i.e.~the number of patterns of the form \vcenteredinclude{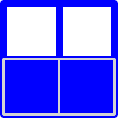}{15pt} (and the rotational analogs) in a configuration;
\item the \emph{number of consecutive pairs}, i.e.~the number of patterns of the form \vcenteredinclude{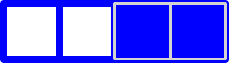}{09pt} (and the vertical analogs) in a configuration.
\end{itemize}
Note that the presence of a very large or very small number of adjacent pairs or consecutive pairs is not expected under the Ising model for any choice of the parameters, and this is an indicator for the presence of positive or negative long-range interactions.

%
%
%

\subsection{Presence of non-homogeneity}
This alternative hypothesis is defined by
%
\begin{equation}\label{E:alternative2}
\prob(y) = \frac{\exp\left(\sum_{i\in V} \alpha_i y_i + \sum_{(i,j)\in E} \beta_{ij}|y_j - y_i| \right)}{Z(\alpha,\beta)},
\end{equation}
where 
$\alpha \in \R^{|V|}$ 
and $\beta \in \R^{|E|}$ are the model parameters and $Z(\alpha, \beta)$ is the normalizing constant. We consider the following test statistics for non-homogeneity: Let  $y\in S$ denote the configuration of the current step of the Markov chain. We uniformly sample $K$ pairs of disjoint square sub-configurations $x_i^1, x_i^2$ of $y$ of size $N\times N$ and compute
%
\begin{equation}\label{E:randsubtables}
{a_i^1} := T_1(x^1), \quad 
{b_i^1}:= T_2(x^1), \quad
{a_i^2} := T_1(x^2), \quad 
{b_i^2}:= T_2(x^2), \quad \textrm{for }\, i=1,\dots ,K.
\end{equation}
The following test statistics measure the degree of non-homogeneity in a configuration:
%
\begin{itemize}
\item \emph{Vertex-non-homogeneity:} $dT_1:=\max_{i=1,\dots ,K} |a_i^1-a_i^2|$. We expect this test statistic to be small under homogeneity when $\alpha_i=\alpha$ for all $i\in V$.

\item \emph{Edge-non-homogeneity:} $dT_2:=\max_{i=1,\dots ,K} |b_i^1-b_i^2|$. We expect this test statistic to be small under homogeneity when $\beta_{ij}=\beta$ for all $(i,j)\in E$.


\item \emph{General non-homogeneity:} $dT_{12} := \max(\frac{\Delta T_1}{N^2}, \frac{\Delta T_2}{2N(N-1)})$, where $N^2$ is the number of vertices and $2N(N{-}1)$ is the number of edges in an $N\times N$ sublattice; this is used for normalization. This test statistic is a combination of the previous two test statistics. 
\end{itemize}
%

\section{Simulations}\label{s:simulations}

In this section, we compare the performance of the test statistics described in Section~\ref{s:TestIsing} for recognizing departure from the Ising model and discuss some simulation results. We performed simulations on the $10\times 10$ lattice and analyzed the type-I and type-II errors of the different test statistics. We generated data under four different models:
\begin{enumerate}
\item[(1)] Ising model
\begin{enumerate}
\item[(i)] with positive interaction, i.e.~$\beta>0$ ,
\item[(ii)] with negative interaction, i.e.~$\beta<0$. 
\end{enumerate}
\item[(2)] Non-Ising model with interaction between diagonal neighbors as described in~\eqref{E:diagonal_interaction}.
\item[(3)] Non-Ising model with vertex-non-homogeneity as described in~\eqref{E:alternative2}.
\item[(4)] Non-Ising model with an overall effect defined by $\alpha=0$, $\beta=0$ and $\gamma=1/5$ as described in~\eqref{E:overall_interaction}.
\end{enumerate}

\begin{figure}[!b]
\centering
\subfigure[]{\includegraphics[scale=0.18]{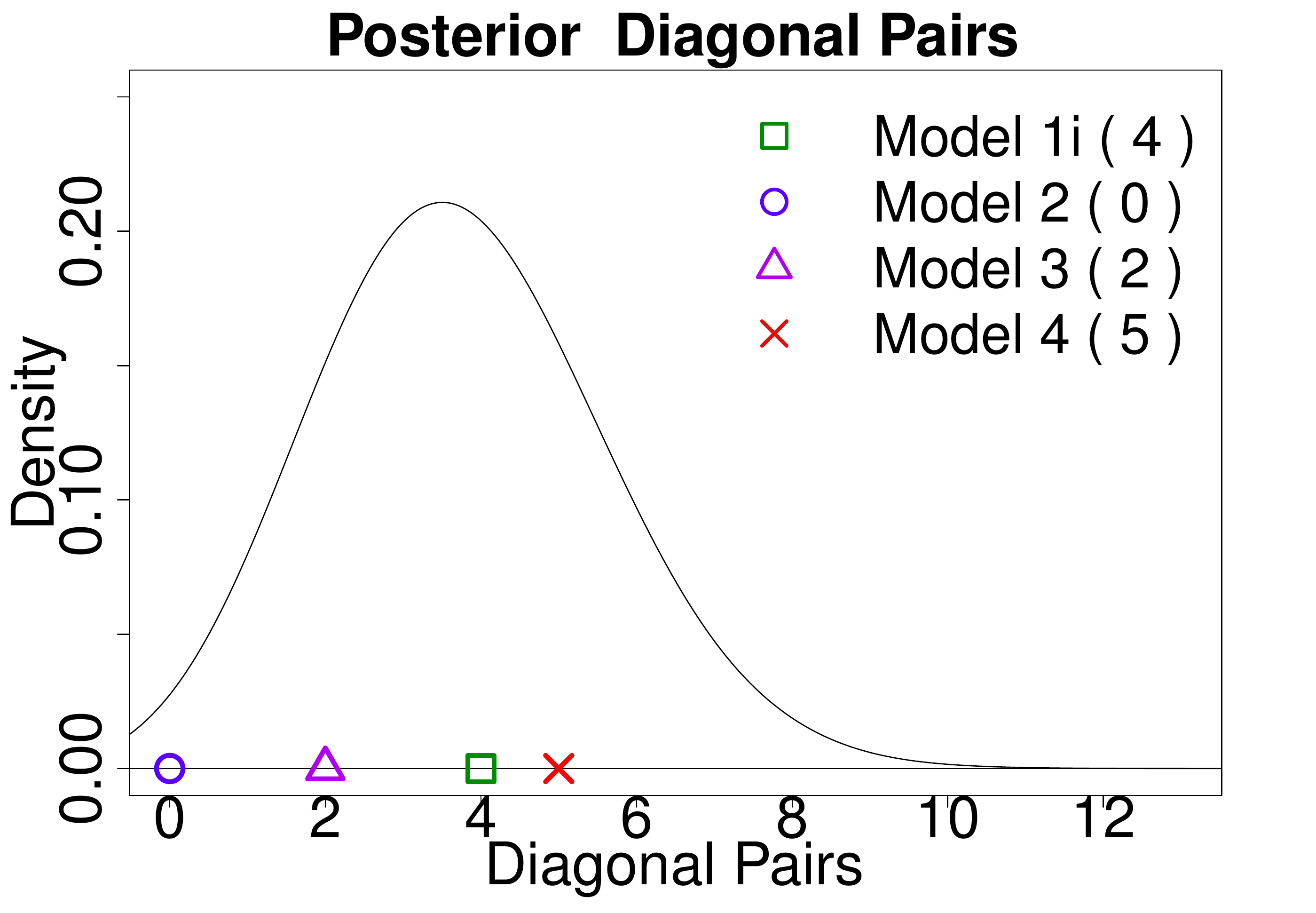}}
\subfigure[]{\includegraphics[scale=0.18]{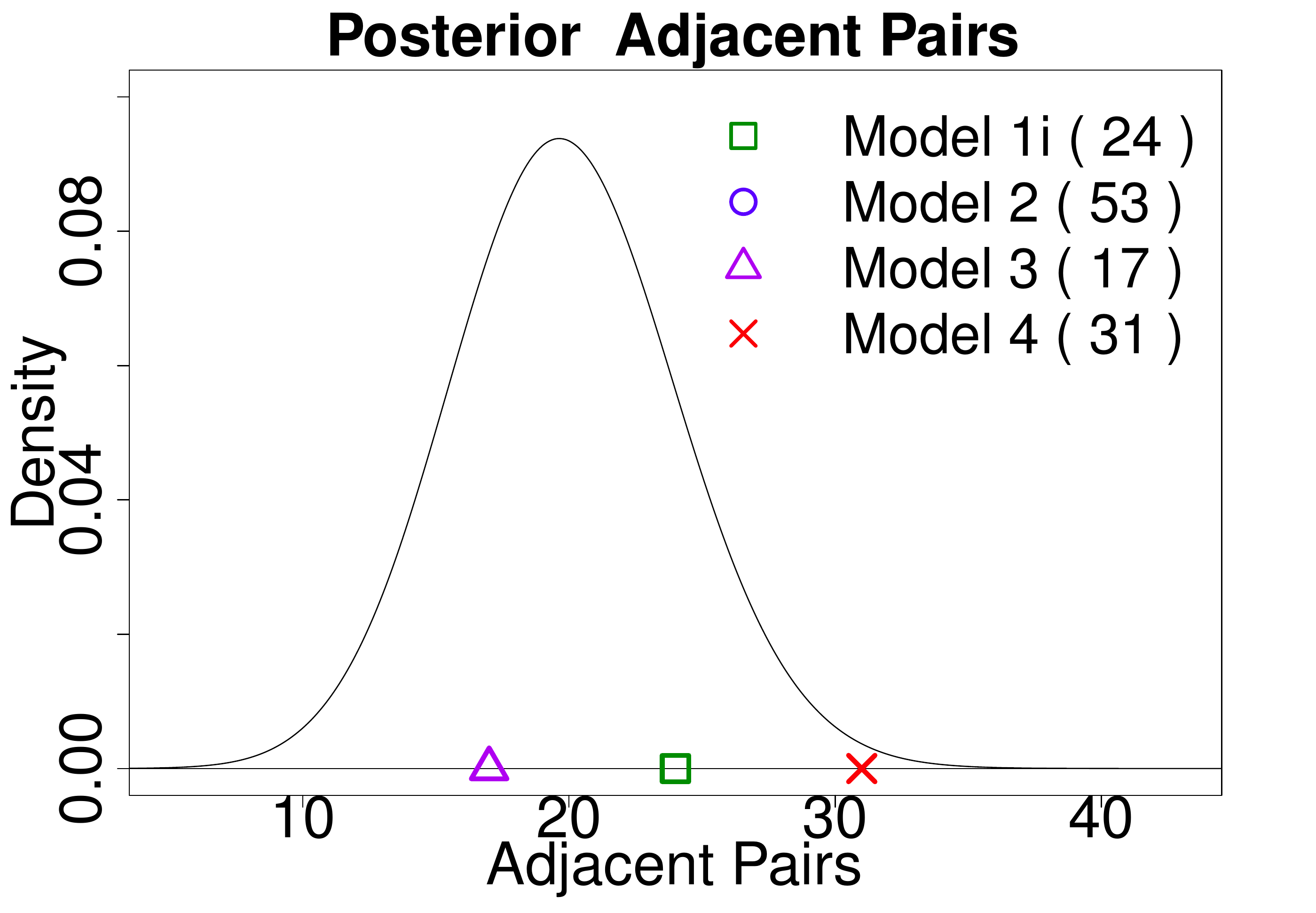}}
\subfigure[]{\includegraphics[scale=0.18]{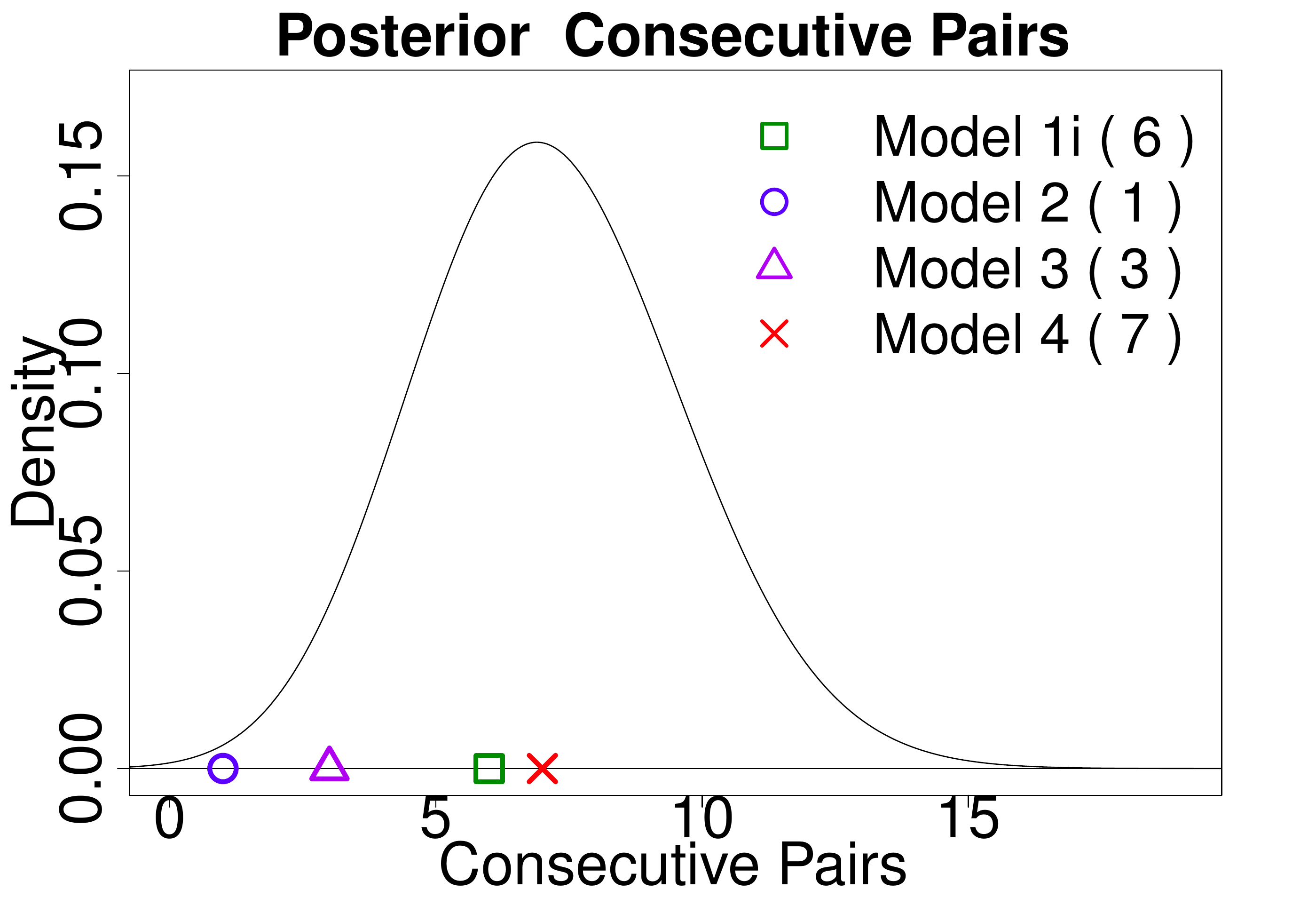}} 
\subfigure[]{\includegraphics[scale=0.18]{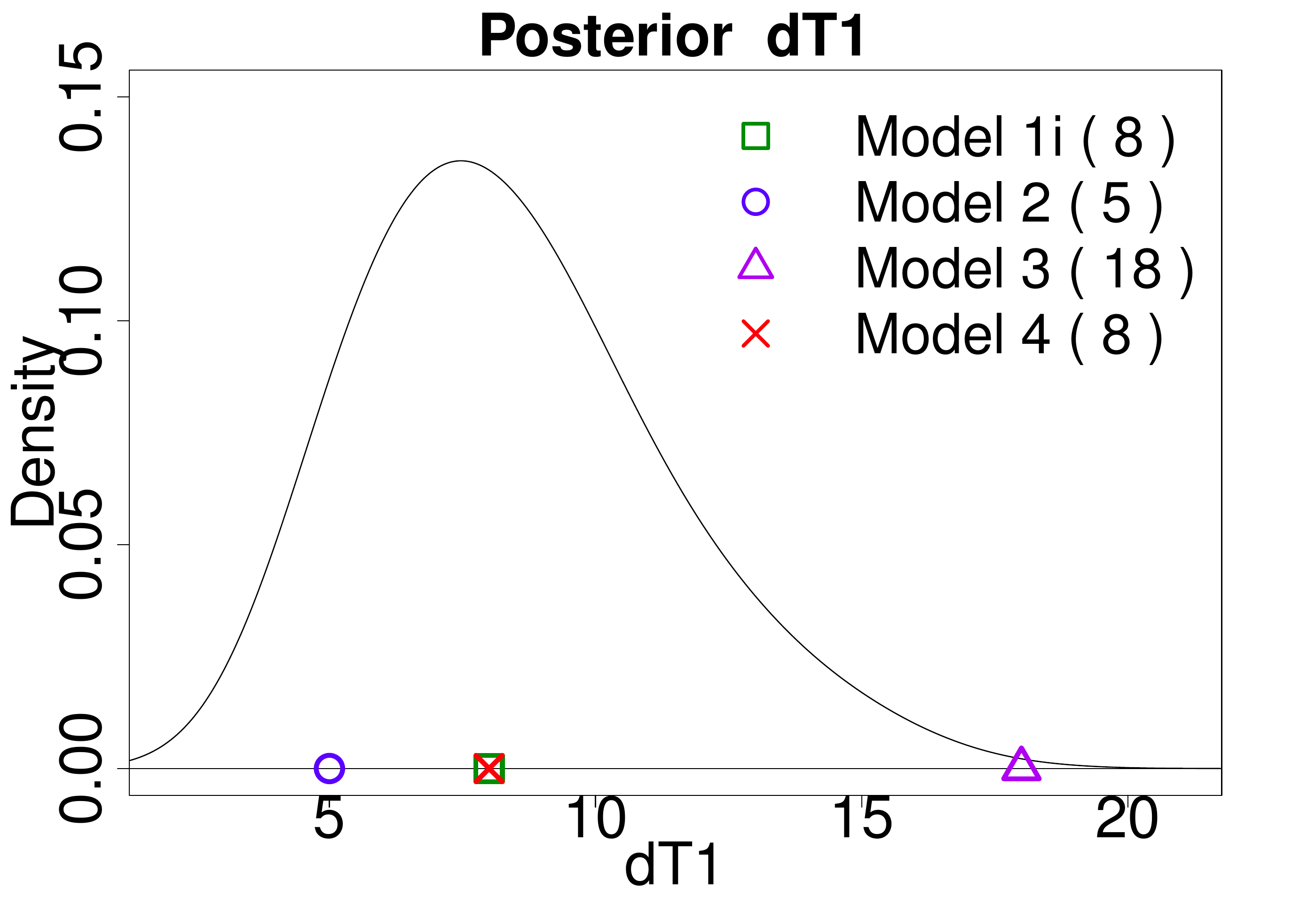}}
\subfigure[]{\includegraphics[scale=0.18]{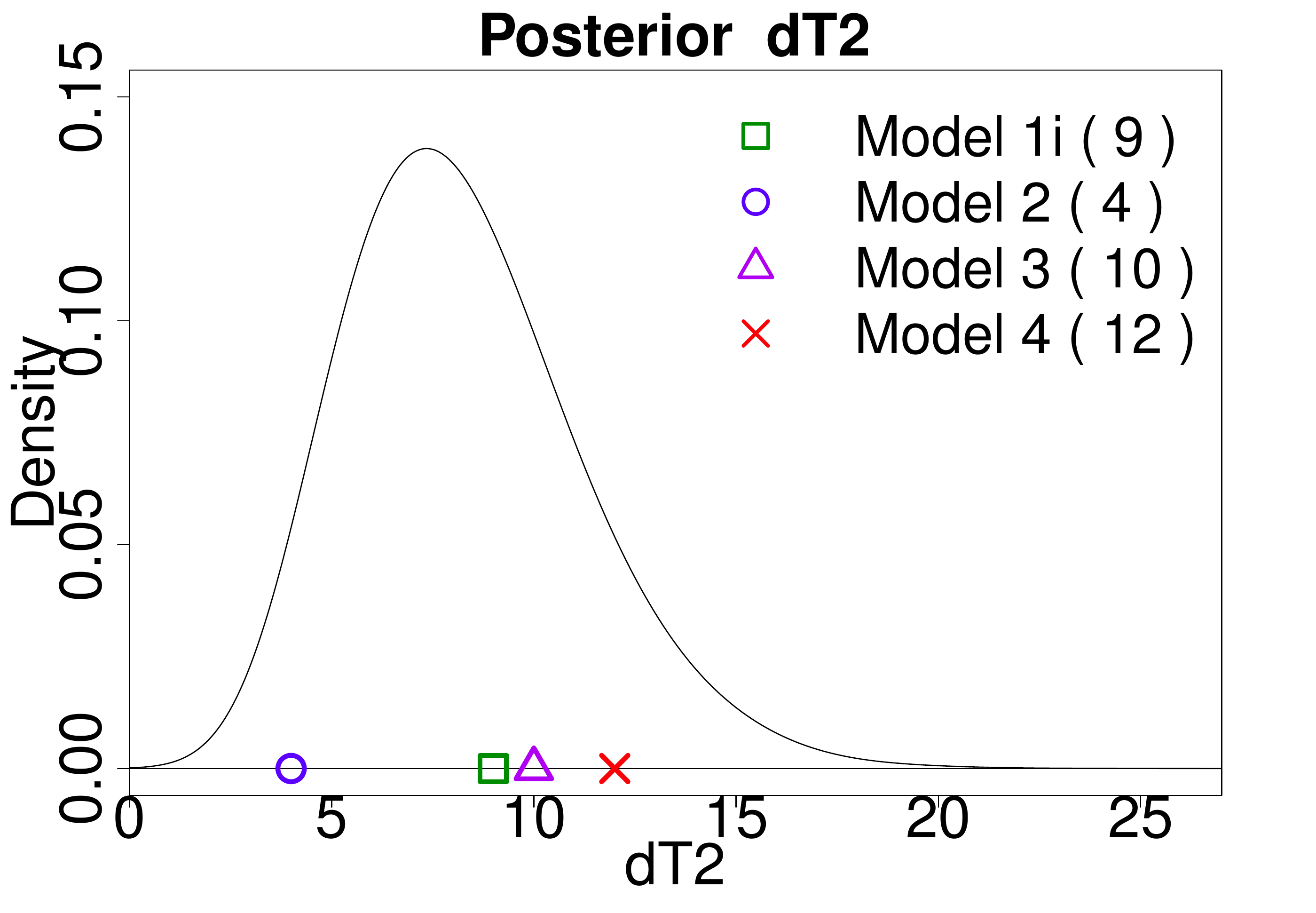}}
\subfigure[]{\includegraphics[scale=0.18]{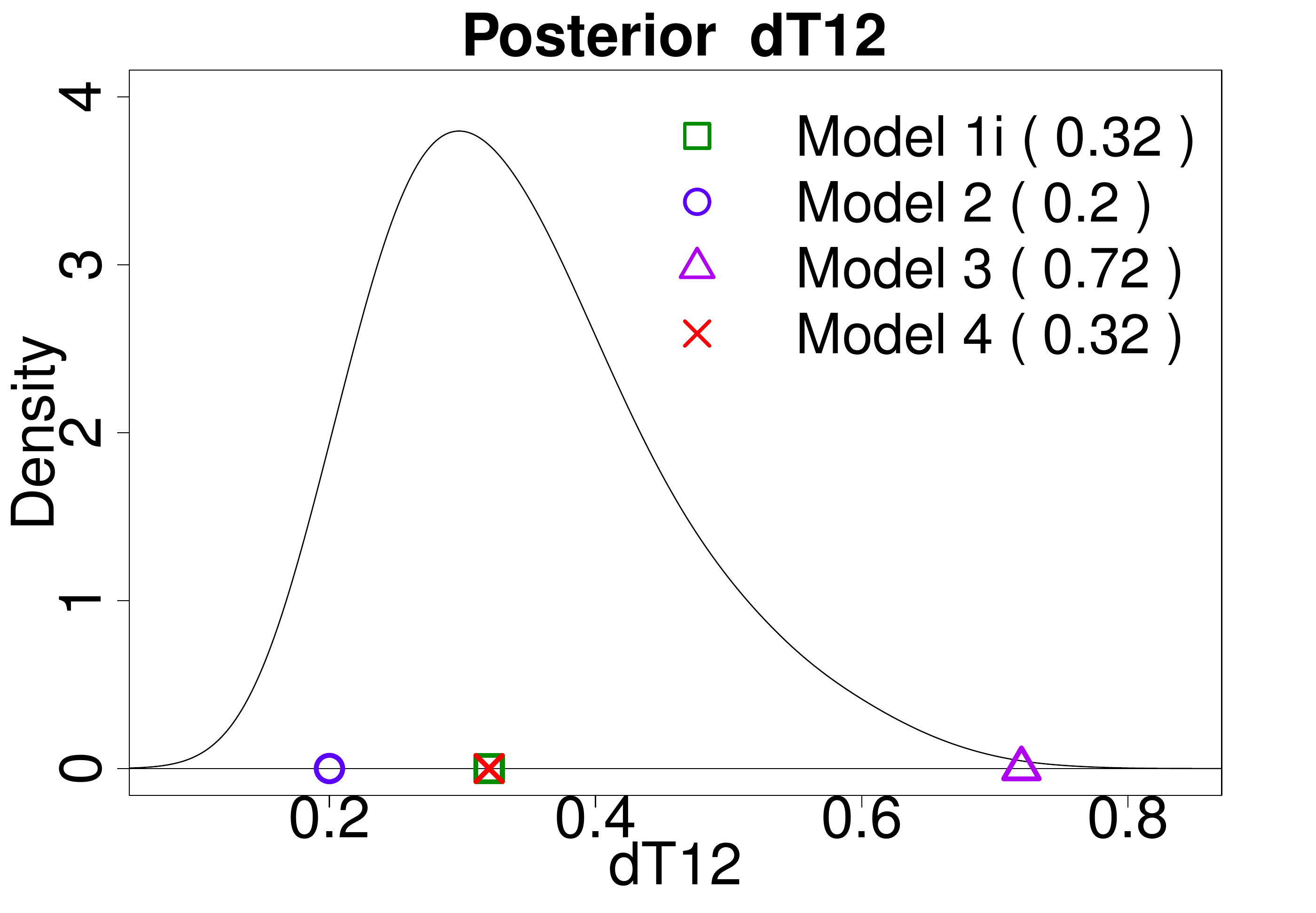}} 
\vspace{-0.3cm}
\caption{Posterior distribution of the six test statistics defined in Section~\ref{s:TestIsing} for the first set of experiments based on three chains of 10,000 Monte Carlo steps each; the sufficient statistics of the observed configurations from the four models 
are $T_1(y) = 52$ and $T_2(y) = 70$; the observed value of each test statistic is given in the legends and depicted in the plots.}\label{F:posterior_set1}
\end{figure}

We generated data from these models by MCMC simulation using periodic boundary conditions. We performed two sets of simulations, namely one set of simulations to compare the test statistics for the Ising model with \emph{positive} interaction to the three non-Ising models and a second set of simulations to compare the test statistics for the Ising model with \emph{negative} interaction to the three non-Ising models. The values for the model parameters were chosen in such a way that, at equilibrium, the resulting data under the different models have the same sufficient statistics; namely $T_1(y)=52$ and $T_2(y)=70$ in the first set of experiments and $T_1(y)=57$ and $T_2(y)=106$ in the second. To generate the data sets, a Markov chain was initiated in a random configuration and the standard Metropolis-Hastings algorithm for minimizing the energy was used for 10,000 steps. After removing the first 1,000 steps as burn-in, we chose at random one of the remaining configurations with the given values of $T_1$ and $T_2$. The data sets for the different models are made available on the website of the first author.

We then generated three Markov chains with 100,000 iterations each as explained in Section~\ref{s:MkvChains} and used the tools described in~\cite{Gilks95} to assess convergence of the chains. This included analyzing the Gelman-Rubin statistic and the autocorrelations. We combined the three Markov chains to generate the posterior distribution of the six test statistics given in Section~\ref{s:TestIsing}. For the non-homogeneity tests, we sampled $K=100$ pairs of subtables of size $3\times 3$.

Figure~\ref{F:posterior_set1} shows the posterior distribution of the six test statistics described in Section~\ref{s:TestIsing} for the first set of simulations with an Ising model with positive interaction and sufficient statistics $T_1(y)=52$ and $T_2(y)=70$. The observed values of the test statistics for the four different models are given in the legends and depicted in the plots. The mean, standard deviation, and various quantiles of the posterior distributions are given in Table~\ref{tab:summary_posteriors_set1}. As seen in Figure~\ref{F:posterior_set1} and in Table~\ref{tab:summary_posteriors_set1} none of the tests reject the null hypothesis for data generated under the Ising model with positive interaction (model 1i) using a significance level of 0.05. For data generated under the non-Ising model with interaction of diagonal neighbors (model 2) the null hypothesis is rejected by all three tests for large-range interaction (a)-(c), but the null hypothesis is not rejected by the tests for non-homogeneity (d)-(f). For data generated under the non-Ising model with vertex-non-homogeneity (model 3) the null hypothesis is rejected using the vertex-non-homogeneity test $dT_1$ and the general non-homogeneity test $dT_{12}$ (Figure~\ref{F:posterior_set1} (d) and (e)) and it is not rejected using any other tests. This shows the importance of non-homogeneity tests. As shown in these simulations, recognizing an overall effect is difficult. Interestingly, test (b) which counts the adjacent pairs is able to recognize departure from the Ising model for this example. 

\begin{table}[b]
\caption{Summary statistics of the posterior distributions for the first set of experiments} 
\centering
\begin{tabular}{c cccccccc}
\hline\hline
Test & Mean & $\sigma$ & \multicolumn{5}{c}{Quantiles} \\
\hline \\
\multicolumn{3}{c}{} & 0.025 & 0.05 &  0.5 & 0.95 & 0.975 \\[0.5ex]
\hline \\
diagonal pairs & 3.727558 & 1.7006851 & 1  & 1  & 4  & 7 &  7\\
adjacent pairs & 19.851085 & 3.9803748 & 13 & 13 & 20 & 27 & 28\\
consecutive pairs & 7.160328 & 2.3170231 &  3  & 4  & 7 & 11 & 12  \\
$dT_1$  &  8.394267 & 2.829628 & 4&  4 &  8 & 14 &   15 \\
$dT_2$  & 8.2569350 & 2.8088039 & 4&  4 &  8 &  13 &   15 \\
$dT_{12}$  &   0.3430878 & 0.1075384 & 0.175 & 0.2 & 0.32 & 0.56  & 0.6 \\  [1ex] 
\hline
\end{tabular}
\label{tab:summary_posteriors_set1}
\end{table} 

\begin{table}[b]
\caption{Summary statistics of the posterior distributions for the second set of experiments} 
\centering
\begin{tabular}{c cccccccc}
\hline\hline
Test & Mean & $\sigma$ & \multicolumn{5}{c}{Quantiles} \\
\hline \\
\multicolumn{3}{c}{} & 0.025 & 0.05 &  0.5 & 0.95 & 0.975 \\[0.5ex]
\hline \\
diagonal pairs & 13.4802400 & 3.4080909 & 7  & 8  & 14  & 19 &  20\\
adjacent pairs & 17.8396450 & 3.8995787 & 10 & 12 & 18 & 25 & 25\\
consecutive pairs & 4.7658575 & 2.2035294 &  1  & 1  & 5 & 9 & 9  \\
$dT_1$  &  5.345140 & 1.703516 & 3&  3 &  5 & 8 &   9 \\
$dT_2$  & 8.9426925 & 3.1193738 & 4&  5 &  9 &  15 & 16 \\
$dT_{12}$  &   0.2543455 & 0.0696624 & 0.15 & 0.16 & 0.24 & 0.375  & 0.4 \\  [1ex] 
\hline
\end{tabular}
\label{tab:summary_posteriors_set2}
\end{table}

\begin{figure}[!t]
\centering
\subfigure[]{\includegraphics[scale=0.18]{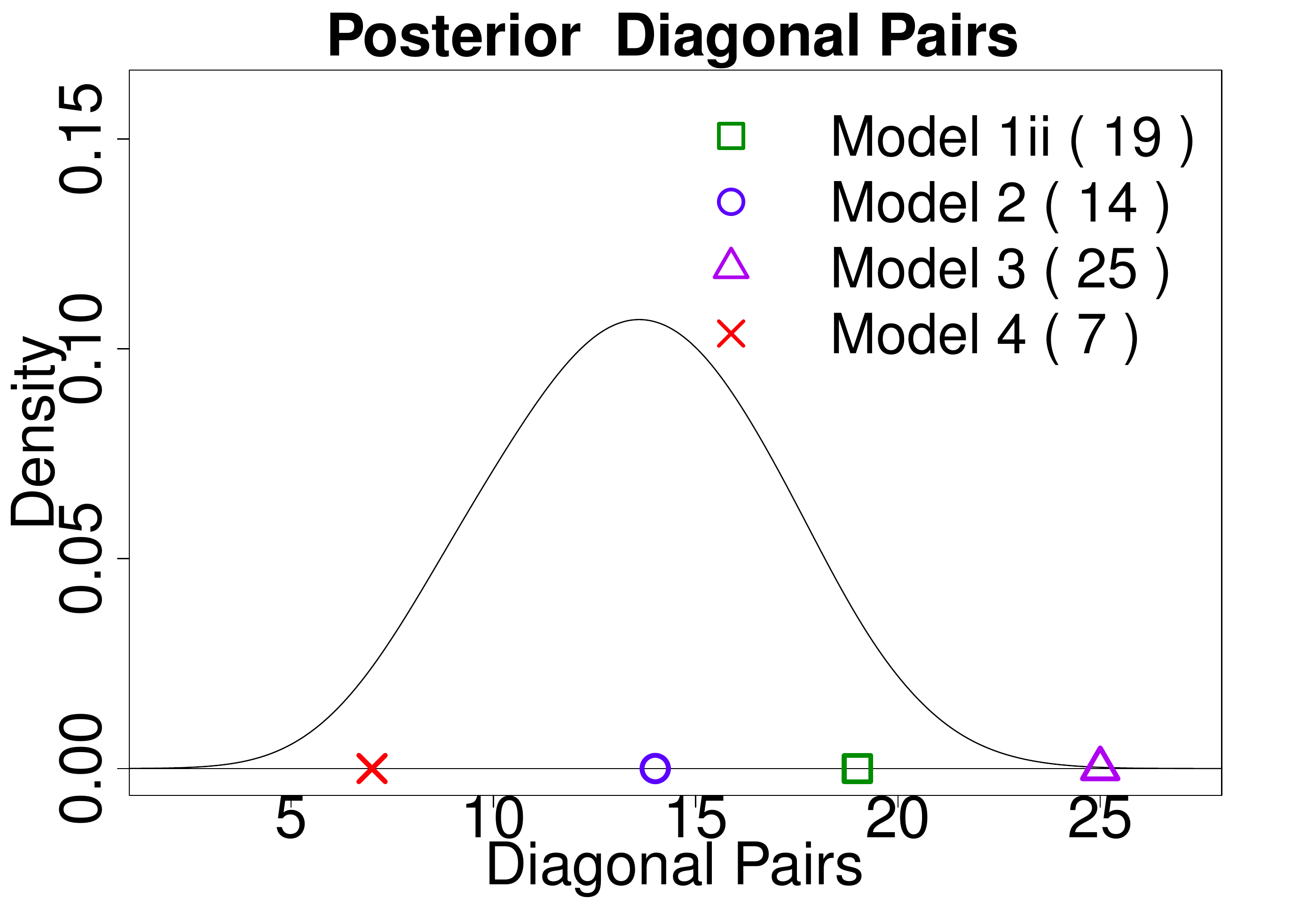}}
\subfigure[]{\includegraphics[scale=0.18]{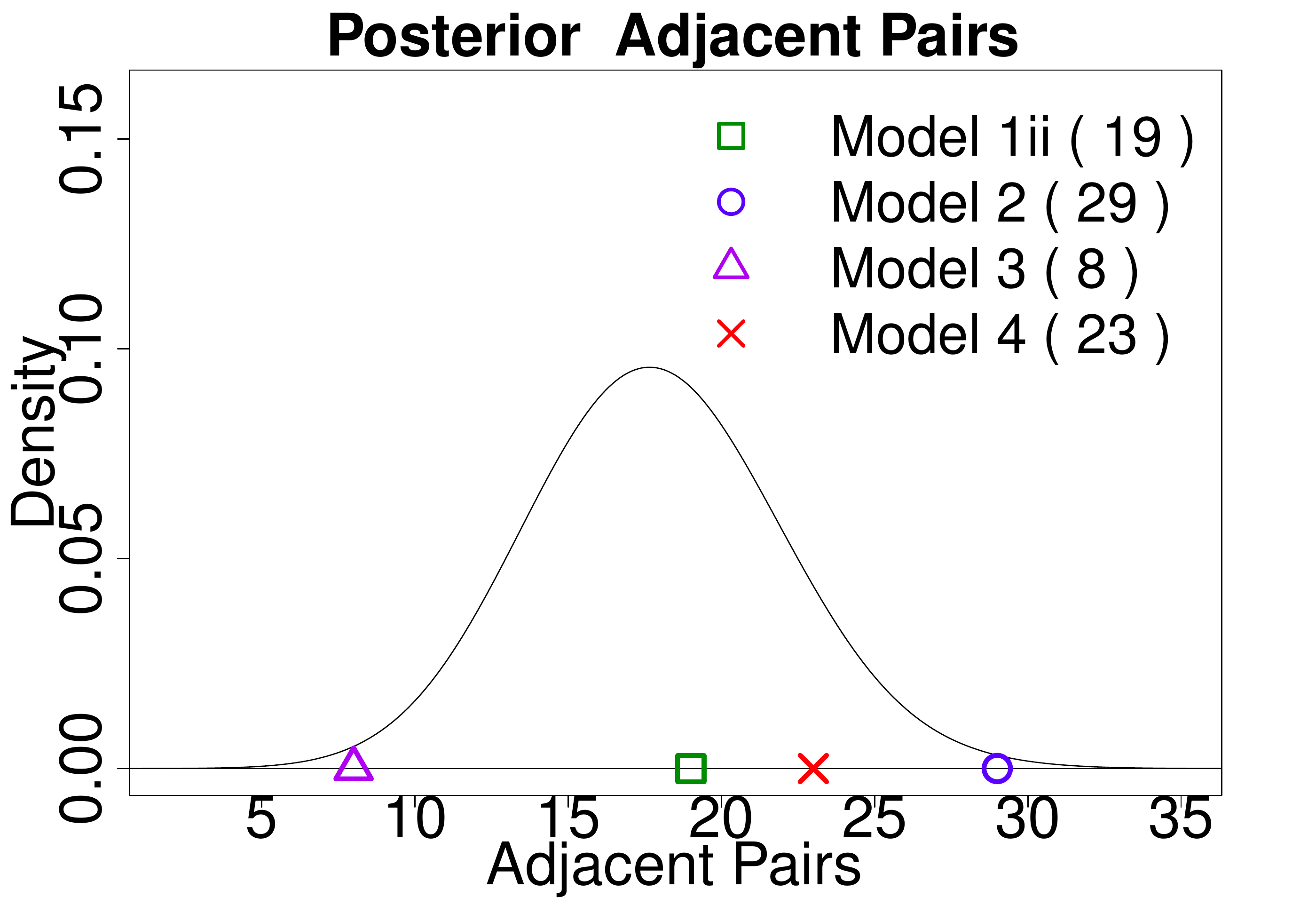}}
\subfigure[]{\includegraphics[scale=0.18]{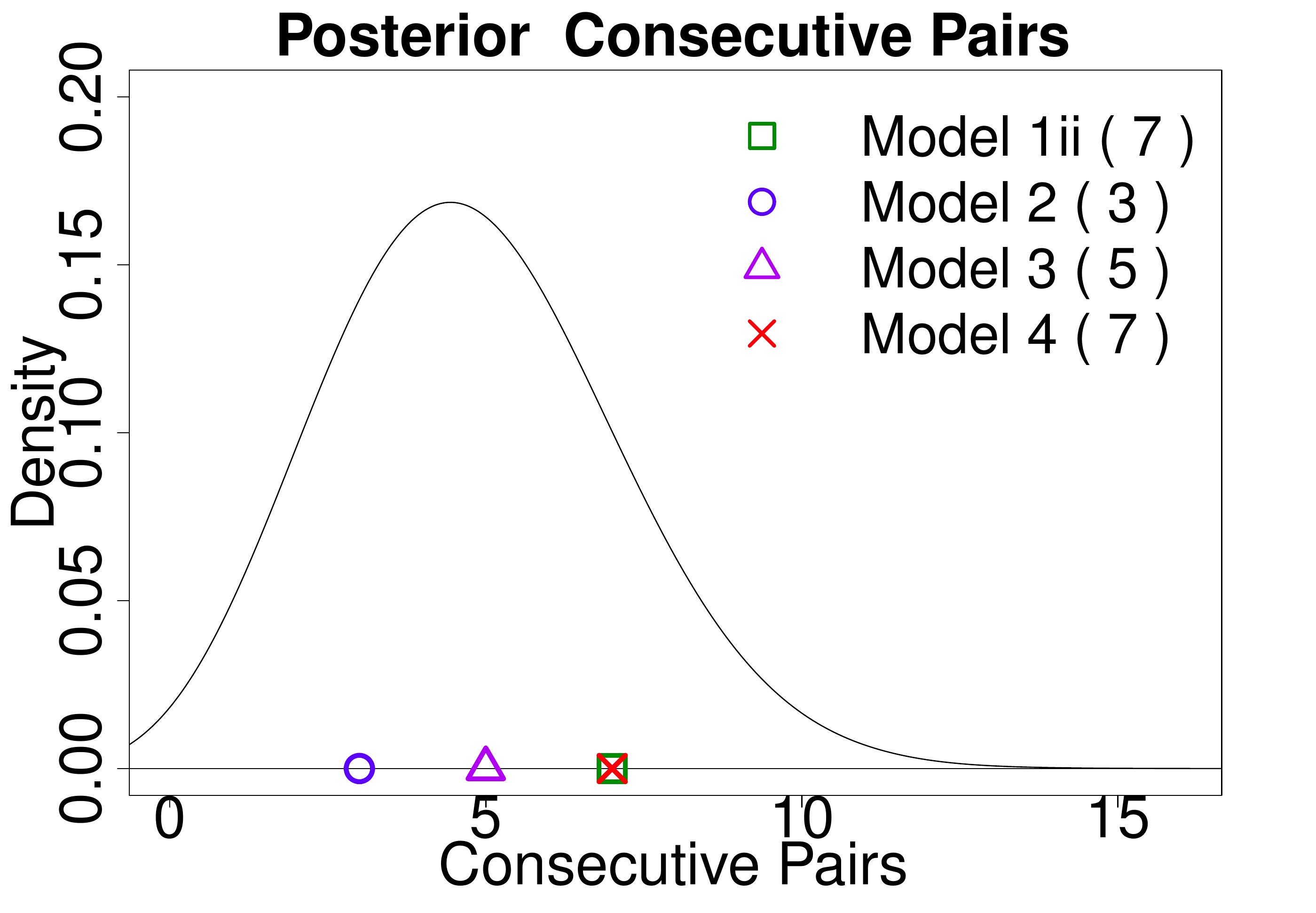}} 
\subfigure[]{\includegraphics[scale=0.18]{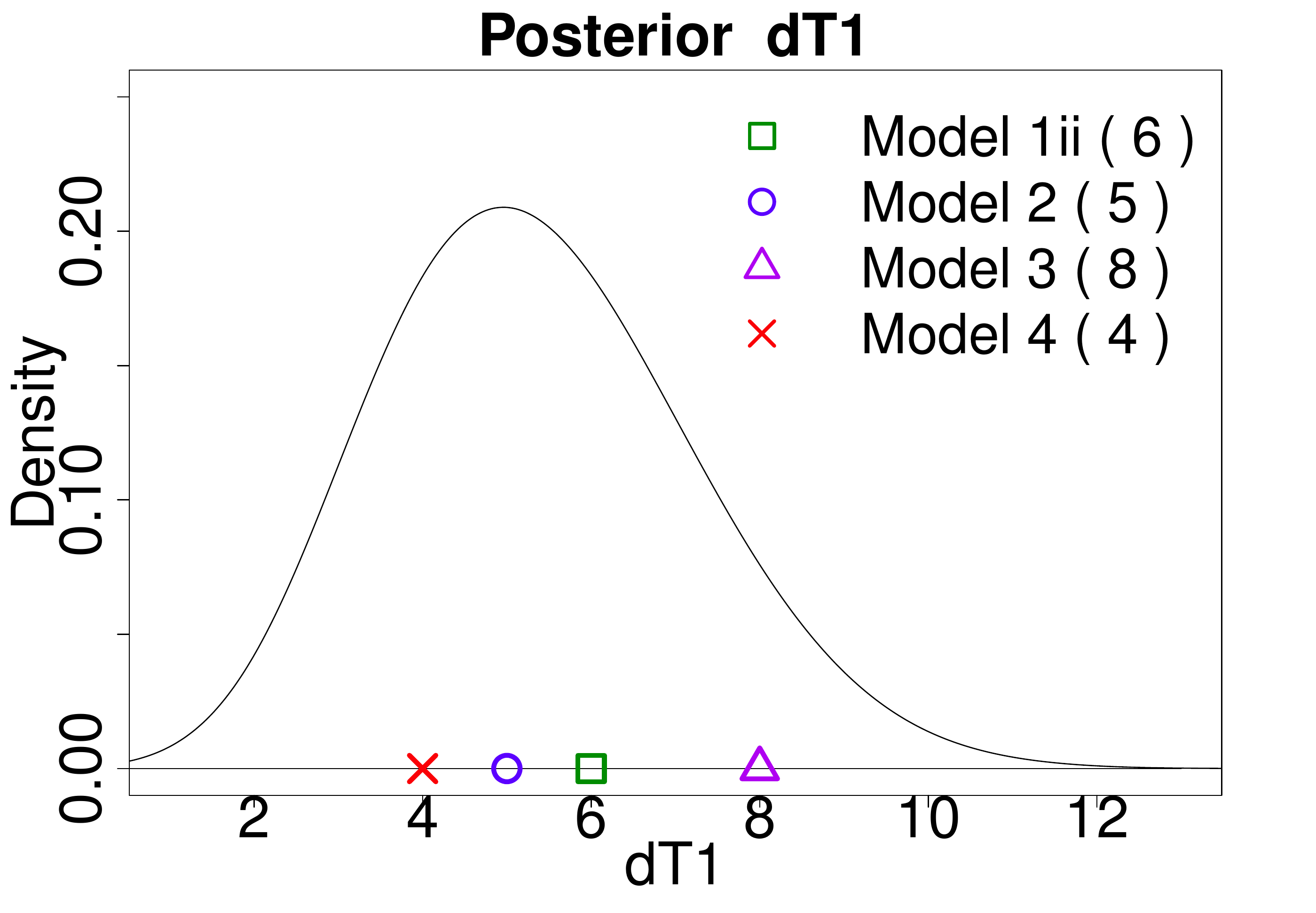}}
\subfigure[]{\includegraphics[scale=0.18]{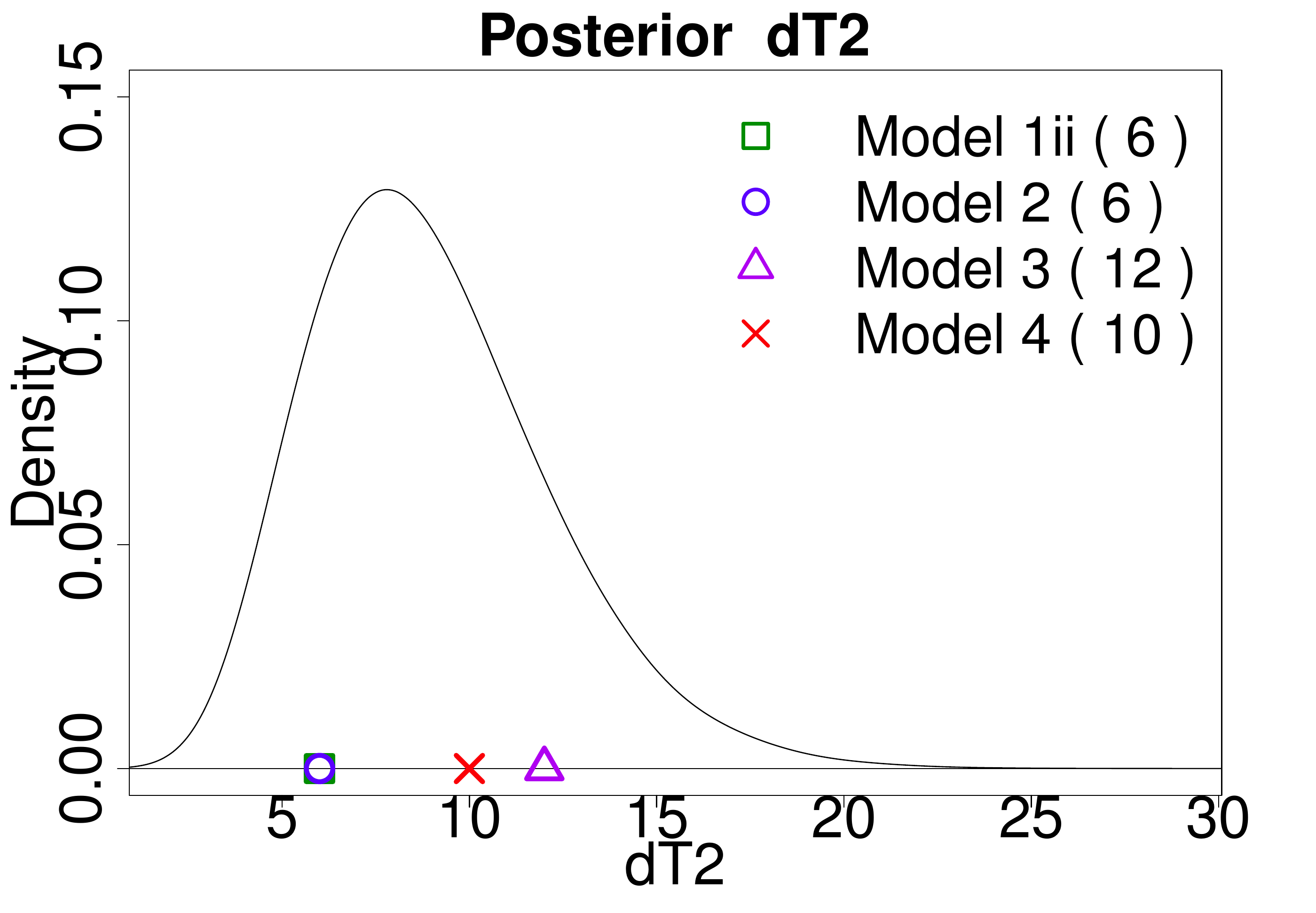}}
\subfigure[]{\includegraphics[scale=0.18]{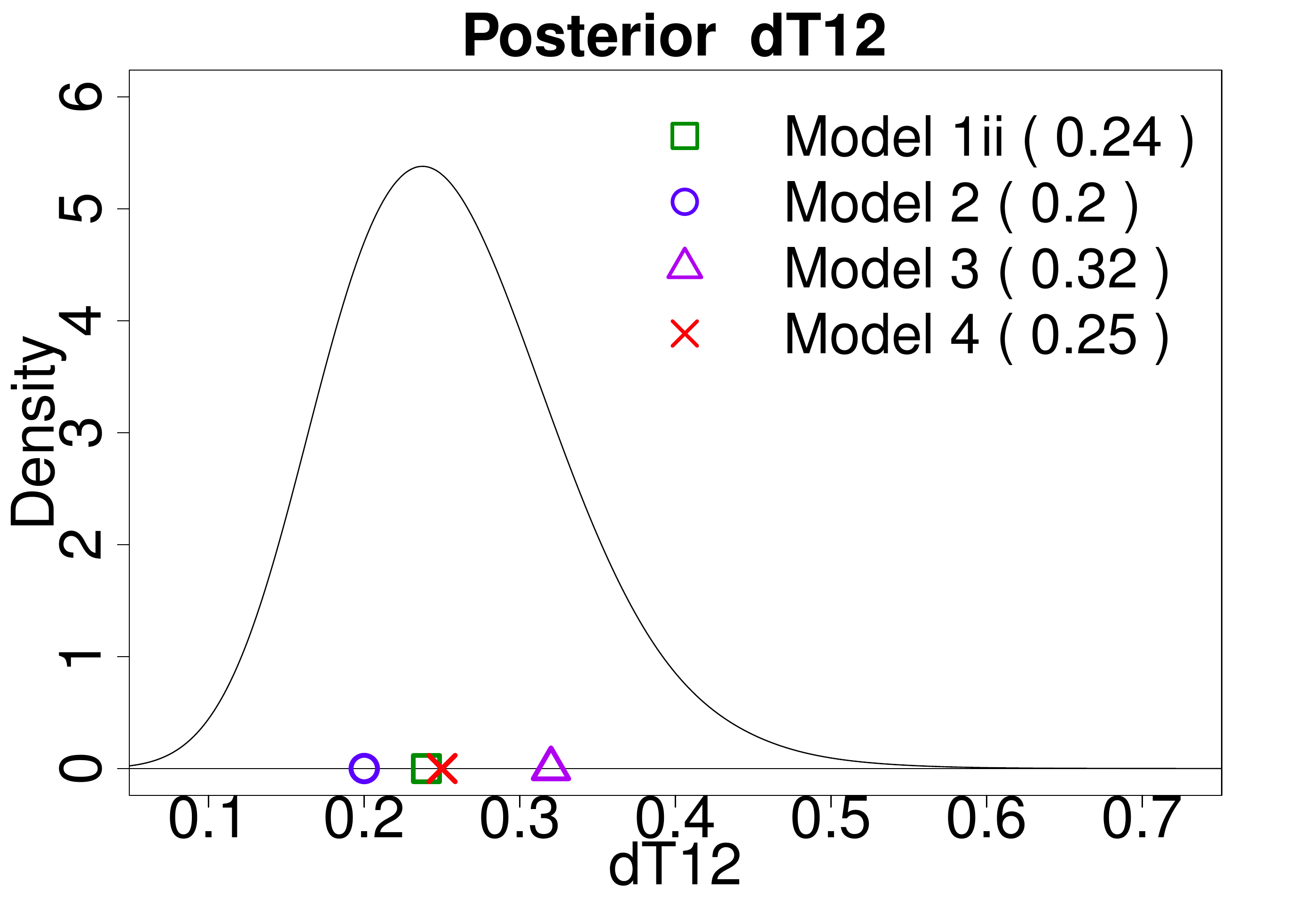}} 
\caption{Posterior distribution of the six test statistics defined in Section~\ref{s:TestIsing} for the second set of experiments based on three chains of 10,000 Monte Carlo steps each; the sufficient statistics of the observed configurations from the four models are $T_1(y) = 57$ and $T_2(y) =106$; the observed value of each test statistic is given in the legends and depicted in the plots.}\label{F:posterior_set2}
\end{figure}

Figure~\ref{F:posterior_set2} shows the posterior distribution of the six test statistics for the second set of simulations with an Ising model with negative interaction and sufficient statistics $T_1(y)=57$ and $T_2(y)=106$. The observed values of the test statistics for the four different models are given in the legends and depicted in the plots. The mean, standard deviation, and various quantiles of the posterior distributions are given in Table~\ref{tab:summary_posteriors_set2}. 
As seen in Figure~\ref{F:posterior_set2} and in Table~\ref{tab:summary_posteriors_set2}, none of the tests reject the null hypothesis for data generated under the Ising model with negative interaction (model 1ii). However, the test based on the number of diagonal pairs is close to rejection even for a moderate negative interaction. For a larger negative interaction this test rejects the null hypothesis leading to a large type-I error. For data generated under the non-Ising model with interaction of diagonal neighbors (model 2) the null hypothesis is only rejected by test~(b) based on counting the number of adjacent pairs. For data generated under the non-Ising model with vertex-non-homogeneity (model 3) the null hypothesis is rejected by two tests for large-range interactions, namely tests (a) and (b) and weakly rejected by the vertex-homogeneity test $dT_1$ (Figure~\ref{F:posterior_set2} (d)). As in the first set of simulations, it is very difficult to recognize an overall effect.

Based on our simulation results we decided to use two test statistics for analyzing the spatial organization of receptors on the cell membrane in Section~\ref{s:biodata}, namely the test statistic based on counting adjacent pairs and the general non-homogeneity test $dT_{12}$. In our simulations, these test statistics seem to have a low type-I and type-II error rate.

\section{Applications to biological data}\label{s:biodata}
In this section, we present an application of our methods to biological data. The data concerns the spatial distribution of receptors on the cell membrane, 
and it consists of an image of a cell membrane in super resolution, where receptors are highlighted against all other components. In order to minimize the border effect that would confer higher density of receptors around the edges in the picture, only a central lattice of $800{\times} 800$ pixels was chosen. This corresponds to the largest square completely enclosed within the circular border. 
To each pixel we associate a  random variable $y_i$ that takes values $0$ or $1$ indicating whether a receptor was present in pixel $i$.  
We apply our goodness-of-fit methodology to test if the spatial distribution of the receptors follows an Ising model.

\begin{figure}[tb]
\centering
\begin{tabular}{rcl}
	\vcenteredinclude{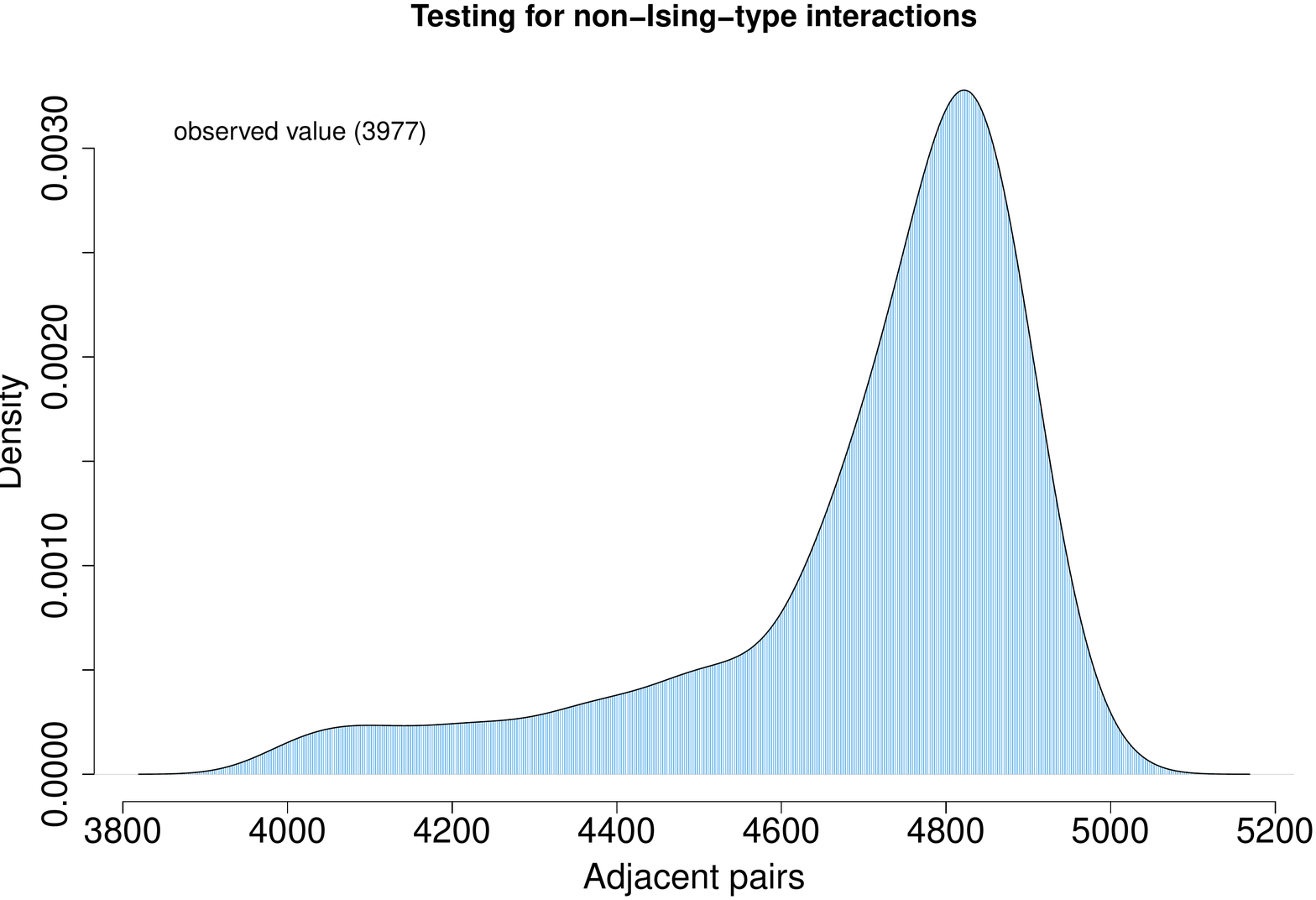}{.35\linewidth}
	\vcenteredinclude{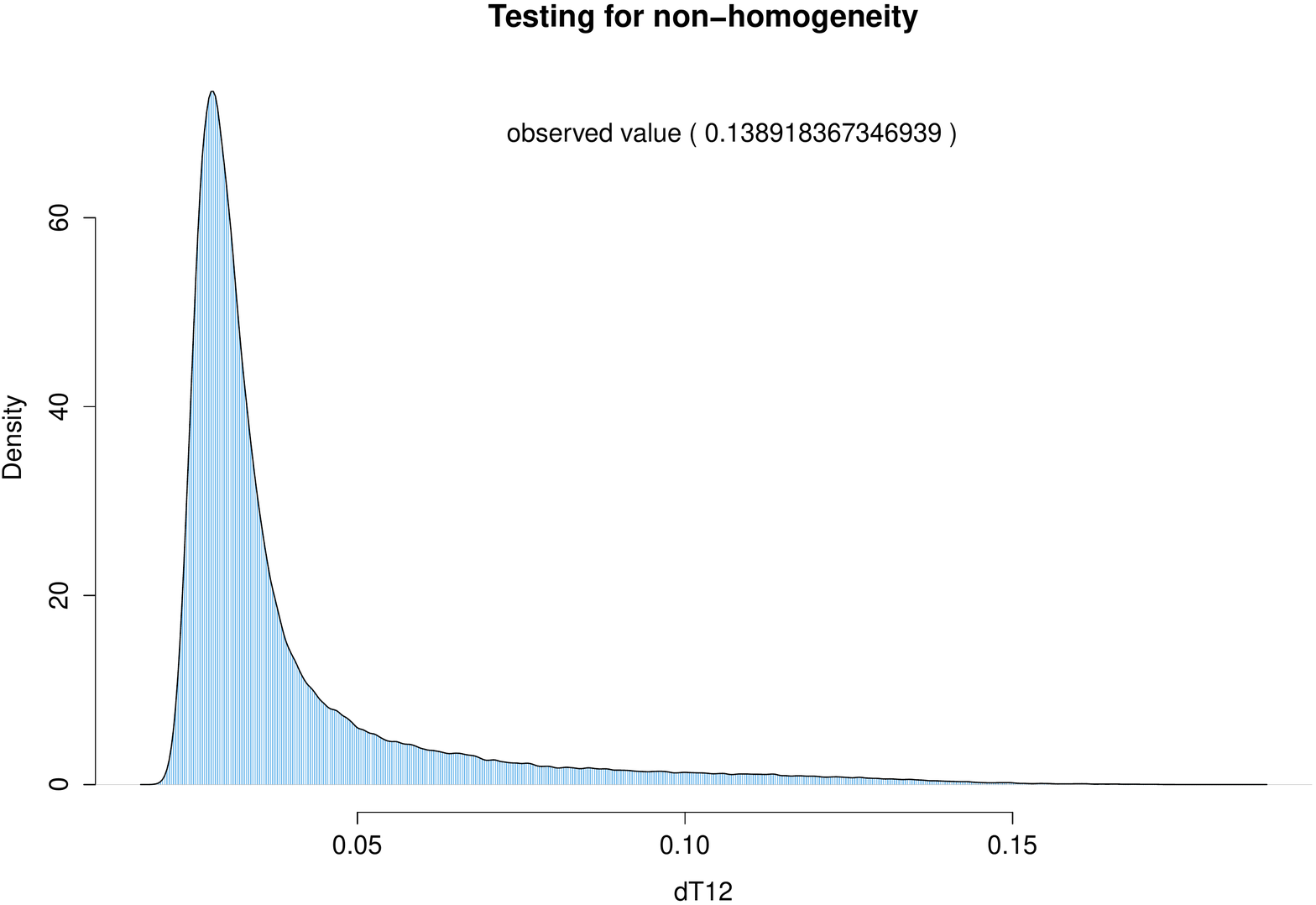}{.37\linewidth}
\end{tabular}
\caption{Posterior distribution of the two test statistics based on 50 chains of 40,000 Monte Carlo steps each starting from the biological data; the observed value of each test statistic is given in the legends.}\label{F:posteriors_cell1}
\end{figure}

The observed configuration $y$ has sufficient statistics $T_1(y) = 14,483$ and $ T_2(y) = 51,145$. As test statistics we used the count of adjacent pairs and the non-homogeneity test $dT_{12}$. The observed values were 3,977 for the adjacent pairs and 0.1389184 for $dT_{12}$.
The Monte Carlo simulation consisted of 50 different realizations of a 40,000 step Markov chain starting from the observed configuration. After removing the first 10,000 steps as burn-in, we analyzed the convergence of the chains by studying the Gelman-Rubin statistic and the autocorrelations as described in~\cite{Gilks95}.
We combined the 50 Markov chains to generate the posterior distribution of the test statistics, illustrated in Figure~\ref{F:posteriors_cell1}. For the non-homogeneity tests, we sampled $K=500$ pairs of subtables of size $50{\times}50$.

The hypothesis of nearest neighbor interaction was discarded at the significance level $\hat p< 0.001$ and the homogeneity hypothesis at significance level $\hat p<0.005$.

%

\appendix
\section{Proofs for the 1-dimensional Ising model}

\subsection{Proof of Lemma~\ref{L:shapeOptimal1D}} 
Suppose that $y\in S(a,b)$ and write  $y=(c,s)$,  where $s$ is the number of singletons of $y$ and $c$ is a configuration without singletons. We need to show that if $y$ is a max-singleton configuration, 
then $c$ consists of at most one connected component of size $|c| = a - b/2 + 1$. 

Assume that $c$ has (at least) two connected components $C_1, C_2$ of size $k_1, k_2 >1$ respectively. Write $C_1=y_{i_1+1}y_{i_1+2} \cdots y_{i_1+k_1}$ and 
$C_2 = y_{i_2+1}y_{i_2+2} \cdots  y_{i_2+k_2}$ with $i_1{+}k_1< i_2$. Without loss of generality, assume that $y_j = 0$ for all $j\in \{i_1{+}k_1{+}1,\ldots i_2-1\}$. Consider the simple swap $z_1={\bf e}_{i_1+k_1+1}-{\bf e}_{i_2}$ and notice that 
$c{+}z$ contains two components, of sizes $k_1{+}1$ and $k_2{-}1$ respectively. Therefore, we can find simple swaps $z_1,\ldots, z_{k_2-1}$ such that $c{+}z_1+\cdots + z_{k_2-1}$ contains one component of size $k_1{+}k_2{-}1$ and one singleton. However, the configuration $y+z_1+\cdots + z_{k_2-1}$ contains one more singleton than $y$, which contradicts the maximality of $y$.
%
%
\hfill \qedsymbol

\subsection{Proof of Proposition~\ref{P:irred1D}}
Write $y=(c,s)$ where $s$ is the number of singletons and $c$ is a configuration without singletons. In particular, we can write 
$c=C_1+\cdots + C_\ell$, where $C_1,\cdots, C_\ell$ are adjacent connected components of length larger than one.
Let $k_i$ denote the length of the component $C_i$ for $i=1,\ldots, \ell$. If $\ell = 0,1$ then $y$ is already an optimal configuration by Lemma~\ref{L:shapeOptimal1D}, thus assume $\ell \geq 2$. As in the proof of Lemma~\ref{L:shapeOptimal1D}, 
we can use simple swaps to decrease the length of $C_\ell$ and increase the length of $C_{\ell-1}$ until $|C_\ell | = 1$ and $|C_{\ell-1}| = k_{\ell-1}+k_\ell -1$. Similarly, we can decrease the length of $C_{\ell-1}$ and increase the length of $C_{\ell-2}$, until the former becomes a singleton and the latter becomes a component of length $k_{\ell-2}+k_{\ell-1}+k_{\ell}-2$. Continuing in this way, we can reduce $C_\ell, C_{\ell-1}, \ldots, C_2$ into singletons and increase $C_1$ to be a component of size $k_{1}+k_{2}+\cdots +k_{\ell}-(\ell -1)$. The resulting configuration is optimal by Lemma~\ref{L:shapeOptimal1D}.
\hfill \qedsymbol

\section{Proofs for the 2-dimensional Ising model}

\subsection{Proof of Lemma~\ref{lemma_opt_2d}}
Let $y^*\in S(a,b)$ be a max-singleton configuration. Let $s$ denote the number of singletons and let $c$ denote the configuration $y^*$ without the singletons. Hence $b = T_2(c)+4s$ and since $s$ is maximal, we are interested in lower bounds on $T_2(c)$. So we solve the following optimization problem:
\begin{equation*}
\begin{split}
&\underset{c}{\mbox{minimize}}  \qquad T_2(c)\\
&\mbox{subject to } \quad\;
T_1(c) = a-s.
\end{split}
\end{equation*}

\begin{figure}[tb]
\begin{tabular}{rcl}
	\vcenteredinclude{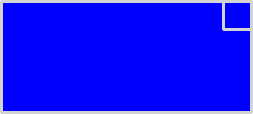}{.09\linewidth}
	& $\longrightarrow$ & 
	\vcenteredinclude{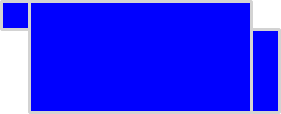}{.09\linewidth}
\end{tabular}
\caption{Turning a configuration $(n,m,0,0,s)$ to $(n{-}1,m,m{-1},1,s)$.}\label{F:breakcorner}
\end{figure}

Let $\mathcal L_c = (\mathcal V_c, \mathcal E_c)$ be the subgraph induced by the configuration $c$. Since $T_2(c) = 4 (a{-}s) - 2| \mathcal E_c|$, minimizing $T_2(c)$ is equivalent to maximizing $|\mathcal E_c|$.  In the language of grid graphs,  Harary and Harborth~\cite{HH76} showed that the maximal number of edges that a graph like $\mathcal L_c$ can have is $\lfloor 2(a{-}s) - 2\sqrt{a{-}s}\rfloor$, and this maximum is attained when $\mathcal L_c$ comes  from a rectangular configuration $c^*$ defined by $(r, m ,d_1,d_2,s)$ with $s=0$, $m = \floor{\sqrt{a-s}}$, $r = \floor{(a-s)/m}$, $d_1=a-s-mr$ and $d_2=0$. As a result, $T_2(c)$ is lower bounded as follows:
\begin{equation}\label{ineq}
4(a-s)-2\lfloor 2(a-s)-2\sqrt{a-s}\rfloor\;\leq\; T_2(c) \;=\; b-4s.
\end{equation}
Since $b$ is always even, we find that (\ref{ineq}) is equivalent to
$$\frac{4a-b}{4} \;\leq\; (a-s)-\sqrt{a-s}.$$
By squaring and solving the quadratic equation we obtain 
$$a-s \;\geq\; \frac{(4a{-}b)/2+1+\sqrt{4a-b+1}}{2}.$$
Since we are interested in the solution that maximizes the number of singletons $s$, we obtain
\begin{equation}\label{eq_sing}
s\; =\; a- \Bigl\lceil \dfrac{(4a-b)/2+1+\sqrt{4a-b+1}}{2}\Bigr\rceil \;=\; \Bigl\lfloor\dfrac{b/2-1-\sqrt{4a-b+1}}{2}\Bigr\rfloor.
\end{equation}
Hence, the max-singleton configuration $y^*$ is a rectangular configuration, where the number of singletons is given in equation (\ref{eq_sing}) and the additional connected component consists of a box of size $m\times r$. We still need to increase the first sufficient statistic by $a-s-mr$ and the second sufficient statistic by $b-4s-2(m{+}r)$. Note that since $s$ is maximized, $b-4s-2(m{+}r)\in\{0,2\}$, thus letting $n=r$ and
$$
(d_1, d_2) =  \left\{ \begin{array}{ll} (a-s-mn, 0) & \textrm{if } b-4s-2(m{+}r)= 0\\ (a-s-mn-1, 1) & \textrm{if } b-4s-2(m{+}r) =2\end{array} \right.,
$$
completes the proof for all cases except for the case $a-s-mr=0$ and $b-4s-2(m{+}r) =2$. For this case, we  reduce the second sufficient statistic as illustrated in Figure~\ref{F:breakcorner}. Thus, the resulting rectangular configuration is of size $(r{-}1, m, m-1, 1,s)$. We complete the proof by letting $n=r-1$ and noting that the condition $a-s-mr=0$ is equivalent to $a{-}s/m$ being integer.
\hfill \qedsymbol

\subsection{Proof of Proposition~\ref{P:rectangle2square}}

\begin{figure}[htb]
	\includegraphics[width=.3\textwidth]{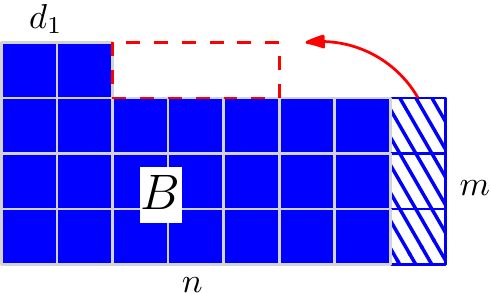}
\vspace{-10pt}
\caption{Changing the shape of the rectangular block $B$.}\label{F:sorting}
\end{figure}

Write $y=(n,m,d_1,d_2,s)$. In addition, suppose that $n-m>1$; otherwise $y$ would already be an optimal rectangular configuration.
We need to show that the $n{\times}m$-rectangular block $B$ can be rearranged into a rectangular block of size $n'{\times}m'$ with $0\leq n'-m'\leq 1$ with only  0- or $\pm2$-swaps.
%
For this, we use moves of the type indicated in Figure~\ref{F:possible moves} to move one of the sides of $B$ of length $m$ to an adjacent side in $B$ as indicated in Figure~\ref{F:sorting}.
\begin{enumerate}
\item If $0 = d_2< d_1$ and $n{-}1{-}d_1 < m$, then let $k = n{-}1{-}d_1$, so we only need $z_1,\ldots, z_k$ moves, all of type A3$\to$B1 (thus they are 0-moves). 
As a result, we get a rectangular configuration $(m',n',d'_1,d'_2,s)$ with
$m' = m{+}1, n'=n{-}1, d'_1=d_1-(n{-}m{-}1)$.

\item If $0=d_2<d_1$ but $n{-}1{-}d_1 \geq m$, then we use moves $z_1,\ldots, z_m$, where $z_1,\ldots, z_{m-1}$ are swaps of type A3$\to$B1 and a move $z_m$ of type A1$\to$B1, which is a $(-2)$-swap. So the resulting configuration will have 
$m'=m, n'=n-1, d'_1 = m+d_1$. Notice that the resulting configuration may not be rectangular if $d_1'>m'$, but we repeat this step until we are in the situation of (1) or (4).

\item If $d_1=0$, 
then 
we just move one of the sides of length $m$ on top, which will require moves $z_1,\ldots, z_m$, where $z_1$ is of type A3$\to$B2,  $z_2,\ldots,z_{m-1}$ of type A3$\to$B1, and $z_m$ of type A1$\to$B1. 

\item If $d_2=1$, then we first make a swap of type A1$\to$B1 to join $d_2$ to the $d_1$ block, giving as a result a configuration with $d'_1 = d_1+1, d_2=0$ but that lie in $S(a,b{-}2)$.
\end{enumerate}
We can always continue in this way, repeating steps (1)--(4), with the subtle difference that if at some point the configuration lies in $S(a,b{-}2)$ 
we will replace the move of type A1$\to$B1 by one of type A1$\to$B3.
In this way, at the end of each step either $n$ decreases by one or $m$ increases by one, resulting in a configuration with $n'-m'\leq 1$. 

If at the end of this procedure the resulting rectangular configuration $ y'=(n,m,d_1,d_2,s')$ lies in $S(a,b{-}2)$, we move to one in $S(a,b)$ as follows. If either $d_2 = 1$ or $d_1=1$, make a simple swap of type A1$\to$B3 to get a configuration in $S(a,b)$ with either $d'_2=0$, or $d'_1=0$ and $s'=s{+}1$. If $d_2=0$ and $d_1\geq 1$, then we make a +2-swap of type A3$\to$B2 to get a configuration with $d'_2 = 1$ and $d'_1 = d_1{-}1$. 

Lastly, if $d_1=0$, then the current configuration has shape $(n,m,0,0,s)$ and with a 2-swap of type A3$\to$B2 we make a configuration of shape $(n{-}1,m,m{-1},1,s)$ as depicted in Figure~\ref{F:breakcorner}.
\hfill \qedsymbol

\bibliographystyle{siam}
\bibliography{ising,revisedbib}

\begin{thebibliography}{10}

\bibitem{4ti2}
{\sc 4ti2 team}, {\em 4ti2---a software package for algebraic, geometric and
  combinatorial problems on linear spaces}.
\newblock {A}vailable at www.4ti2.de.

\bibitem{AHT12}
{\sc S.~S. Aoki, H.~Hara, and A.~Takemura}, {\em Markov {B}ases in {A}lgebraic
  {S}tatistics}, vol.~199 of Springer Series in Statistics, Springer, 2012.

\bibitem{Bes72}
{\sc J.~Besag}, {\em Nearest-neighbour systems and the auto-logistic model for
  binary data}, J. Roy. Statist. Soc. Ser. B, 34 (1972), pp.~75--83.

\bibitem{Bes74}
\leavevmode\vrule height 2pt depth -1.6pt width 23pt, {\em Spatial interaction
  and the statistical analysis of lattice systems}, J. Roy. Statist. Soc. Ser.
  B, 36 (1974), pp.~192--236.
\newblock With discussion by D. R. Cox, A. G. Hawkes, P. Clifford, P. Whittle,
  K. Ord, R. Mead, J. M. Hammersley, and M. S. Bartlett and with a reply by the
  author.

\bibitem{Bes77}
\leavevmode\vrule height 2pt depth -1.6pt width 23pt, {\em Some methods of
  statistical analysis for spatial data}, in Proceedings of the 41st {S}ession
  of the {I}nternational {S}tatistical {I}nstitute ({N}ew {D}elhi, 1977),
  {V}ol. 2, vol.~47, 1977, pp.~77--91, 138--147.
\newblock With discussion.

\bibitem{Besag_image1}
\leavevmode\vrule height 2pt depth -1.6pt width 23pt, {\em On the statistical
  analysis of dirty pictures}, J. Roy. Statist. Soc. Ser. B, 48 (1986),
  pp.~259--302.

\bibitem{Besag_image2}
\leavevmode\vrule height 2pt depth -1.6pt width 23pt, {\em Towards {B}ayesian
  image analysis}, J. Appl. Statist., 20 (1993), pp.~107--119.

\bibitem{BC89}
{\sc J.~Besag and P.~Clifford}, {\em Generalized {M}onte {C}arlo significance
  tests}, Biometrika, 76 (1989), pp.~633--642.

\bibitem{Ising_history}
{\sc S.~G. Brush}, {\em History of the {L}enz-{I}sing model}, Rev. Mod. Phys.,
  39 (1967), pp.~883--893.

\bibitem{BB00}
{\sc F.~Bunea and J.~Besag}, {\em M{CMC} in {$I\times J\times K$} contingency
  tables}, in Monte {C}arlo methods ({T}oronto, {ON}, 1998), vol.~26 of Fields
  Inst. Commun., Amer. Math. Soc., Providence, RI, 2000, pp.~25--36.

\bibitem{cai2014}
{\sc S.~Cai, B.~Li, and J.~Guo}, {\em A simplification of computing {M}arkov
  bases for graphical models whose underlying graphs are suspensions of
  graphs}, Stat. Sinica, 24 (2014), pp.~447--461.

\bibitem{CDY10}
{\sc Y.~Chen, I.~H. Dinwoodie, and R.~Yoshida}, {\em Markov chains, quotient
  ideals and connectivity with positive margins}, in Algebraic and Geometric
  Methods in Statistics, Cambridge Univ. Press, Cambridge, 2010, pp.~99--110.

\bibitem{DS98}
{\sc P.~Diaconis and B.~Sturmfels}, {\em Algebraic algorithms for sampling from
  conditional distributions}, Ann. Statist., 26 (1998), pp.~363--397.

\bibitem{Dob12}
{\sc A.~Dobra}, {\em Dynamic {M}arkov bases}, J. Comput. Graph. Statist., 21
  (2012), pp.~496--517.

\bibitem{DSS09}
{\sc M.~Drton, B.~Sturmfels, and S.~Sullivant}, {\em Lectures on {A}lgebraic
  {S}tatistics}, vol.~39 of Oberwolfach Seminars, Springer, 2009.

\bibitem{Galam_1}
{\sc S.~Galam, Y.~Gefen, and Y.~Shapir}, {\em Sociophysics: {A} new approach of
  sociological collective behaviour. {I}. mean--behaviour description of a
  strike}, J. Math. Sociol., 9 (1982), pp.~1--13.

\bibitem{Galam_2}
{\sc S.~Galam and S.~Moscovici}, {\em Towards a theory of collective phenomena:
  {C}onsensus and attitude changes in groups}, Eur. J. Soc. Psychol., 21
  (1991), pp.~49--74.

\bibitem{Gilks95}
{\sc W.~R. Gilks, S.~Richardson, and D.~J. Spiegelhalter}, eds., {\em Markov
  {C}hain {M}onte {C}arlo in {P}ractice}, Interdiscip. Statist., Chapman \&
  Hall, London, 1996.

\bibitem{GPS14}
{\sc E.~Gross, S.~Petrovi{\'c}, and D.~Stasi}, {\em Goodness-of-fit for
  log-linear network models: {D}ynamic {M}arkov bases using hypergraphs}.
\newblock To appear in Ann.~Inst.~Stat.~Math., preprint arXiv:1401.4896, 2014.

\bibitem{HH76}
{\sc F.~Harary and H.~Harborth}, {\em Extremal animals}, J. Combin. Inform.
  System Sci., 1 (1976), pp.~1--8.

\bibitem{HoSu07}
{\sc S.~Ho{\c{s}}ten and S.~Sullivant}, {\em A finiteness theorem for {M}arkov
  bases of hierarchical models}, J. Comb. Theory (Series A), 114 (2007),
  pp.~311--321.

\bibitem{Ising}
{\sc E.~Ising}, {\em Beitrag zur {T}heorie des {F}erromagnetismus}, Z. Phy., 31
  (1925), pp.~253--258.

\bibitem{MBD}
{\sc T.~Kahle and J.~Rauh}, {\em The {M}arkov bases database}, 2013.
\newblock {\tt http://markov-bases.de/}.

\bibitem{KOT15}
{\sc T.~{Koyama}, M.~{Ogawa}, and A.~{Takemura}}, {\em Markov degree of
  configurations defined by fibers of a configuration}, J. Algebr. Stat., 6
  (2015), pp.~80--107.

\bibitem{Majewski2001}
{\sc J.~Majewski, H.~Li, and J.~Ott}, {\em The {I}sing model in physics and
  statistical genetics}, Am. J. Hum. Genet., 69 (2001), pp.~853--862.

\bibitem{OHT13}
{\sc M.~Ogawa, H.~Hara, and A.~Takemura}, {\em Graver basis for an undirected
  graph and its application to testing the beta model of random graphs}, Ann.
  Inst. Statist. Math., 65 (2013), pp.~191--212.

\bibitem{PRF10}
{\sc S.~Petrovi{\'c}, A.~Rinaldo, and S.~E. Fienberg}, {\em Algebraic
  statistics for a directed random graph model with reciprocation}, in
  Algebraic Methods in Statistics and Probability {II}, vol.~516 of Contemp.
  Math., Amer. Math. Soc., Providence, RI, 2010, pp.~261--283.

\bibitem{PP13}
{\sc G.~Pistone and M.~Rogantin}, {\em The algebra of reversible {M}arkov
  chains}, Ann. Inst. Stat. Math., 65 (2013), pp.~269--293.

\bibitem{PR12}
{\sc G.~Pistone and M.~P. Rogantin}, {\em Toric statistical models: {I}sing and
  {M}arkov}, in Harmony of {G}r\"obner bases and the modern industrial society,
  World Sci. Publ., Hackensack, NJ, 2012, pp.~288--313.

\bibitem{Rap03}
{\sc F.~Rapallo}, {\em Algebraic {M}arkov bases and {MCMC} for two-way
  contingency tables}, Scand. J. Stat., 30 (2003), pp.~385--397.

\bibitem{RY10}
{\sc F.~Rapallo and R.~Yoshida}, {\em Markov bases and subbases for bounded
  contingency tables}, Ann. Inst. Statist. Math., 62 (2010), pp.~785--805.

\bibitem{RS16}
{\sc J.~Rauh and S.~Sullivant}, {\em Lifting {M}arkov bases and higher
  codimension toric fiber products}, J. Symbolic Comput., 74 (2016),
  pp.~276--307.

\bibitem{Schneidman2006}
{\sc E.~Schneidman, M.~J. Berry, R.~Segev, and W.~Bialek}, {\em Weak pairwise
  correlations imply strongly correlated network states in a neural
  population}, Nature, 440 (2006), pp.~1007--1012.

\bibitem{SZP14}
{\sc A.~B. {Slavkovi{\'c}}, X.~{Zhu}, and S.~{Petrovi{\'c}}}, {\em Fibers of
  multi-way contingency tables given conditionals: relation to marginals, cell
  bounds and {M}arkov bases}, Ann. Inst. Stat. Math,  (2014), pp.~1--28.

\bibitem{TMASBB14}
{\sc G.~Tkacik, O.~Marre, D.~Amodei, E.~Schneidman, W.~Bialek, and M.~J. Berry,
  II}, {\em Searching for collective behavior in a large network of sensory
  neurons}, PLoS Comput.~Biol., 10 (2014), p.~e1003408.

\bibitem{Tkacik}
{\sc G.~Tkacik, E.~Schneidman, M.~J. {Berry~II}, and W.~Bialek}, {\em Ising
  models for networks of real neurons}.
\newblock Preprint arXiv:0611072, 2006.

\end{thebibliography}

\end{document}